\documentclass[a4paper, 11pt]{amsart}
\usepackage[latin1]{inputenc}
\usepackage[ T1]{fontenc}
\usepackage[english]{babel}
\usepackage{amssymb}
\usepackage{amsmath}
\usepackage{amsthm}
\usepackage{amscd}
\usepackage{amsfonts}
\usepackage{stmaryrd}
\usepackage{pb-diagram}
\usepackage{epic,eepic,epsfig}
\usepackage{a4wide}
\usepackage{nextpage}
\usepackage{fancyhdr}
\usepackage{enumerate}
\usepackage{hyperref}
\usepackage{float}
\usepackage{comment}
\usepackage{todo}
\usepackage{faktor}
\usepackage{listings}
\usepackage[
backend=biber,
style=alphabetic,
]{biblatex}

\newtheorem{teo}{Theorem}[section]
\newtheorem{defi}[teo]{Definition}
\newtheorem{lema}[teo]{Lemma}
\newtheorem{prop}[teo]{Proposition}
\newtheorem{cor}[teo]{Corollary}
\newtheorem{rem}[teo]{Remark}
\newtheorem{ex}[teo]{Example}

\newtheorem{conj}[teo]{Conjecture}
\newtheorem{quest}{Question}
\newtheorem{vacio}{Assumption}

\renewcommand{\Im}{\operatorname{Im}}
\newcommand{\R}{\mathbb{R}}
\newcommand{\Z}{\mathbb{Z}}
\newcommand{\Q}{\mathbb{Q}}
\newcommand{\Qbar}{\overline{\Q}}

\newcommand{\N}{\mathbb{Z}_{> 0}}
\newcommand{\C}{\mathbb{C}}
\newcommand{\A}{\mathcal{A}}

\newcommand{\eps}{\varepsilon}

\newcommand{\End}{\operatorname{End}}
\renewcommand{\O}{\mathcal{O}}
\newcommand{\al}{\alpha}
\newcommand{\be}{\beta}

\renewcommand{\H}{\mathbb{H}}
\newcommand{\M}{\mathcal{M}}

\newcommand{\SL}{\operatorname{SL}}

\newcommand{\la}{\lambda}
\newcommand{\Hum}{\mathcal{H}}
\renewcommand{\t}{\tau}
\newcommand{\bt}{\boldsymbol{\tau}}
\newcommand{\bz}{\boldsymbol{z}}
\newcommand{\trdeg}{\operatorname{trdeg}}
\newcommand{\rank}{\operatorname{rank}}
\newcommand{\spn}{\operatorname{span}}
\newcommand{\D}{\mathbb{D}}
\newcommand{\Sp}{\operatorname{Sp}}
\newcommand{\spbt}{\spn_\Q (1, \tau_1, \tau_2, \tau_3, \t_2^2 - \t_1 \t_3)}
\newcommand{\trqjq}{\trdeg \Q (q_1, q_2, q_3, j_1(q), j_2(q), j_3(q))}
\newcommand{\trqjqtil}{\trdeg \Q (\tilde{q}_1, \tilde{q}_2, \tilde{q}_3, j_1(\tilde{q}), j_2(\tilde{q}), j_3(\tilde{q}))}
\renewcommand{\L}{\mathcal{L}}
\newcommand{\Llin}{\mathcal{L}^{lin}}
\newcommand{\tr}{\operatorname{tr}}
\newcommand{\Disc}{\operatorname{Disc}}
\newcommand{\Hom}{\operatorname{Hom}}
\newcommand{\NS}{\operatorname{NS}}
\newcommand{\Jac}{\operatorname{Jac}}
\renewcommand{\Re}{\operatorname{Re}}
\newcommand{\No}{\operatorname{N}}
\newcommand{\nrd}{\operatorname{nrd}}
\newcommand{\GL}{\operatorname{GL}}
\newcommand{\Pic}{\operatorname{Pic}}
\begin{document}

\title{On the CM exception to a generalization of the St\'ephanois theorem}
\author{Desir\'ee Gij\'on G\'omez}

\address{Department of Mathematical Sciences, University of Copenhagen,
Universitetsparken 5, 
2100 Copenhagen O, Denmark.}
\email{dgg@math.ku.dk}

\thispagestyle{empty}

\begin{abstract}
There are two classical theorems related to algebraic values of the j-invariant: Schneider's theorem and the St\'ephanois theorem. Schneider's theorem for the j-invariant states that the transcendence degree $\trdeg \Q(\tau, j(\tau)) \geq 1$ with the sole exception of CM points. In contrast, CM points do not constitute an exception to the St\'ephanois theorem, which states $\trdeg \Q(q,j(q))\geq 1$ for the Fourier expansion ($q$-expansion) of the j-invariant, for any $q$. Schneider's theorem has been generalized to higher dimensions, and in particular holds for the Igusa invariants of a genus 2 curve. These functions have Fourier expansions, but a result of St\'ephanois type is unknown. 
In this paper, we find that there are positive dimensional sources of exceptions to the generic behaviour expected in genus 2, and we discuss their relation to CM points. We utilize Humbert singular relations, putting them into the transcendental framework. The computations of the transcendence degree for CM points are conditional to Schanuel's conjecture.
\end{abstract}

 \maketitle
{\flushleft
}

\begin{center}
---------
\end{center}

\tableofcontents

\section{Introduction}

Consider the elliptic j-invariant $j: \H \to \C$, where $\H = \{ \tau \in \C \lvert \; \Im \t >0 \}$, a $\SL_2(\Z)$-invariant modular function, which classifies elliptic curves over $\Qbar$. It admits a Fourier expansion:
\[
j(q) = \frac{1}{q} + 744 + 196884q + \cdots \in \Z[[q]], \, q = e^{2\pi i \tau},
\]
which defines a meromorphic function on $\D = \{ q \in \C:$ $|q| < 1\}$.

For a field $\Q \subset L$ we denote $\trdeg L$ for the transcendence degree of $L$ over $\Q.$ The following two results are classical. 

\begin{teo}[Schneider's theorem \cite{Schneider}] For $\tau \in \H$, $\trdeg \Q (\tau, j(\tau)) \geq 1$ with the sole exception of $\tau$ quadratic imaginary.
\end{teo}

\begin{teo}[The St\'ephanois theorem\footnote{This is the standard name for this result in French, but to the best of our knowledge this theorem does not have a consistent name in English, save for the previous denomination as the Mahler-Manin conjecture. The French name comes from the result being proven in St-\'Etienne and "st\'ephanois" is its gentilic in French.} \cite{Stephanois}] For $0 < |q| < 1$, $\trdeg\Q(q, j(q)) \geq 1$.
\end{teo}

The $j$-invariant admits a generalization for curves of genus two, which we call Igusa invariants $j_1, j_2, j_3$, see \cite{IguAVM}. Remark that, in the literature, Igusa invariants, or Igusa-Clebsch invariants, usually refer to weighted projective invariants, while we mean absolute ones. Also note that in \cite{IguAVM} there are five weighted invariants instead of four, so that the theory extends to fields of even characteristic.

A curve of genus two is necessarily hyperelliptic and admits a (singular) model $C :y^2 = f(x)$ for $f$ a polynomial of degree $6$. Set $\al_i$ for its six complex roots. We define the following:

\begin{align*}
I_2(f) &= \sum (\al_1 - \al_2)^2 (\al_3 - \al_4)^2 (\al_5 - \al_6)^2, \\
I_4(f) &= \sum (\al_1 - \al_2)^2 (\al_2 - \al_3)^2 (\al_3 - \al_1)^2 (\al_4 - \al_5)^2 (\al_5 - \al_6)^2 (\al_6 - \al_4)^2, \\
I_6(f) &= \sum \prod_{l=1}^3\prod_{(i,j) \in C_l}(\al_i - \al_j)^2, \\
\text{where } & C_1 = \{(1,2),(2,3),(3,1)\}, C_2=\{(4,5),(5,6),(6,4)\}, C_3 = \{(1,4),(2,5),(3,6)\},\\
I_{10}(f) &= \prod (\al_i - \al_j)^2,
\end{align*}
where the sums range among all permutations of six elements. Note that $I_{10}$ is the discriminant of $f$. Therefore, as we are considering \emph{smooth} curves of genus two, we assume in the following that $I_{10} \not = 0$.

We set the following normalization of the Igusa invariants
\begin{align*}
j_1(C) &= \frac{I_2^5}{I_{10}}, \\
j_2(C) &= \frac{I_2^3 I_4}{I_{10}}, \\
j_3(C) &= \frac{I_2^2 I_6}{I_{10}},
\end{align*}
if $I_2 \not= 0$. There are other normalizations to choose when $I_2 = 0$, but we would not need them in this paper. As $I_{\al}$ are symmetric polynomial functions on the roots of $f$, they are functions on the coefficients of $f$. This mimics the situation with the $j$-invariant of an elliptic curve, which also admits a formula in terms of its Weierstrass model. They characterize completely the geometric isomorphism class of the curve (\cite[Section 6, Theorem 2]{IguAVM}), and they detect algebraicity of genus two curves: a curve $C$ is defined over $\Qbar$ if and only if $j_l(C) \in \Qbar$ for $l=1,2,3.$ In contrast to the case of elliptic curves, this does not hold over a fixed number field $K$: if $j_i(C) \in K$ for $i=1,2,3$, then $C$ may not be defined over $K$, but over a quadratic extension of $K$ (see \cite{MestreConstructionDeCourbes}). 

However, we need to consider the Igusa invariants as analytic functions on the moduli space of genus two curves $\mathcal{M}_2$. These analytic functions are not directly defined on $\mathcal{M}_{2}$, but on its image under the Torelli map in $\mathcal{A}_2$ (Torelli locus), the moduli space of principally polarized abelian surfaces. In other words, they can also be seen as invariants of the Jacobian of the curves (endowed the principal polarization). By the Torelli theorem (see \cite[Theorem 12.1]{SilvCornelArithGeom}), over any algebraically closed field this map is injective, or equivalently, the Jacobian variety equipped with the principal polarization completely distinguishes the curve.

It is a classical result (\cite[Satz 2]{Weil}, or \cite[Corollary 11.8.2 a)]{Birk-LangeCAV}) that the Torelli locus is precisely the indecomposable locus of $\A_2$, which we denote $\A_2^{ind} = \A_2 \smallsetminus (\A_1 \times \A_1)$, where for $\A_1 \times \A_1$ we mean the identification with the locus of products of elliptic curves with the product principal polarization. We remark that an abelian surface in the indecomposable locus can be isogenous to a product of elliptic curves, or even isomorphic, if considering another principal polarization in the product of elliptic curves that is not induced as a product polarization.

The moduli space $\A_2$ is coarsely represented by a quotient of a symmetric space, analogously to $\A_1$ with respect to $\SL_2(\Z) \backslash \H$. The analytic space is the Siegel upper half-space of degree two:
\[
\H_2 = \{ \bt \in \M_2(\C)| \: \bt = \bt^{t}, \: \Im \bt \text{ is positive definite} \}.
\]
In the following, we set $\bt =  \begin{pmatrix}
\tau_1 & \tau_2 \\ \tau_2 & \tau_3
\end{pmatrix}$. Note that $\Im \bt$ is positive definite if and only if:
\[
\tau_1, \tau_3 \in \H, \: \Im(\tau_2)^2 < \Im(\tau_1)\Im(\tau_3). 
\]
We denote the symplectic group
\[
\Sp_4(\Z) = \left\{ M \in \M_4(\Z)|\;  M^t J M = J, \text{ where } J =\left( \begin{array} {c | c}
0 & I_2  \\ \hline  -I_2 & 0
\end{array} \right) \right\},
\]
which acts on $\H_2$ by linear fractional transformations: if $M = \begin{pmatrix}
    \al & \be \\ \gamma & \delta
\end{pmatrix} \in \Sp_4(\Z)$ then 
\[
M\bt = (\al \bt + \be)(\gamma \bt + \delta)^{-1}.
\]

The theory of modular forms generalizes to these functions, which are Siegel modular functions. In particular, they admit higher dimensional Fourier expansions (as they satisfy analogous relations to $j(\tau + 1) = j(\tau)$) of the form
\[
\sum_{\substack{n \in \operatorname{GL}_2(\Q) \\ n \text{ half integral symmetric}}} a(n) e^{2\pi i \tr(n\bt)},
\]
where $\tr(n\bt)= n_1 \tau_1 + n_3 \tau_3 + 2n_2 \tau_2$. Therefore, they admit series expansions in terms of $q_l = q_l(\bt) = e^{2\pi i \tau_l},$ for $l=1,2,3.$ 

As we are interested in geometric invariants of abelian surfaces, they are always considered as defined over an algebraically closed field, and all our related definitions (morphisms, isogenies, simplicity, etc.) are geometric as well. 

The analog of Schneider's theorem is a known result for every genus, and in more general settings.

\begin{teo}[Cohen-Shiga-Wolfart\footnote{The set-up of abelian \emph{surfaces} with quaternionic multiplication was historically first study in \cite{Morita}.} \cite{ShigaWolfart} and \cite{CohenHumSurfTransc}] Consider a principally polarized abelian variety $A$ with a (small) period matrix $\bt$. Then the following are equivalent:
\begin{itemize}
\item the abelian variety $A$ is defined over $\Qbar$ and $\bt \in \M_{g}(\Qbar)$,
\item the abelian variety $A$ has complex multiplication.
\end{itemize} \label{TeoCohenShigaWolfart}
\end{teo} 

In our set-up, it implies that the components of the matrix $\bt$ and $j_l(\bt)$ for $l=1,2,3$ are all simultaneously algebraic if and only if $\bt$ parametrizes a CM abelian surface.

We recall the definition of a CM (complex multiplication) abelian surface; they are the natural generalization of CM elliptic curves. By Poincar\'e's irreducibility theorem, an abelian variety $A$ admits an isogeny to a product
\[
\prod_{k=1}^{s} A_k^{r_k} 
\]
for $A_k$ simple abelian varieties. We use the notation $A \simeq B$ for isogenous abelian varieties. We denote $\End(A)$ for the (geometric) endomorphism ring and $\End_0(A) = \End(A) \otimes \Q$ for the endomorphism algebra.

\begin{defi} We say that a \emph{simple} abelian variety $A$ has CM if its endomorphism algebra $\End_0(A) = K$ for $K$ a CM field with $[K : \Q] = 2\dim (A)$. An abelian variety is said to have CM if all the simple factors given by its isogeny decomposition have CM. We call $\bt \in \H_g$ a CM \emph{point} if it parametrizes a CM abelian variety.
\end{defi}

We also remark that this property does not involve the polarization of the abelian variety.

An abelian surface $A$ with CM is necessarily one of the following. We follow the notation from \cite[Table 1]{Daw-OrrUI}.
\begin{itemize}
\item Either $A$ is simple, then $\End_0(A)$ is a quartic CM field \emph{(simple)};
\item or $A \simeq E \times E'$ for two not isogenous $CM$ elliptic curves (with necessarily different CM fields $K,K'$), then $\End_0(A) = K \times K'$ \emph{(non-simple)};
\item or $A \simeq E^2$, then $\End_0(A) = \mathcal{M}_2(K)$ \emph{(isotypic)}.
\end{itemize}

We note that in this notation, a CM point that is not simple can be \emph{non-simple} or \emph{isotypic.}

\subsection{Questions and our main result}
The analog of the St\'ephanois theorem is still open, but the corresponding functional transcendence statement does hold: it follows from \cite[Theorem 2]{BerZud}. (On that direction of functional transcendence, see \cite[Theorem 2.5]{PilaModALWderivtives}.) One could tentatively predict a "generic" behaviour
\begin{equation}
\trqjq \geq 3, \label{EqGenBehav}
\end{equation}
where $q = (q_1,q_2,q_3)$, that would be satisfied outside a set of exceptions. 

The purpose of this paper is to answer two questions, posed to us by Daniel Bertrand.
\begin{quest} Can CM points be exceptions to \eqref{EqGenBehav}?
\end{quest}

\begin{quest} In case of an affirmative answer, are such exceptions exclusively due to being CM, or are they explained by belonging to a larger "exceptional subset"?
\end{quest}

This last question makes sense from the perspective of Shimura varieties, and the work behind the proof of the Andr\'e-Oort conjecture. Likewise, the functional transcendence statement in \cite[Theorem 2.5]{PilaModALWderivtives} has weakly special subvarieties as unique sources of exceptions, as well as \cite[Theorem 1.1] {CGFNAxLindWderivatives} and \cite[Theorem 1.3]{PilTsimAxSforjfunction}. Also note that, under the Grothendieck's period conjecture, in \cite[Proposition 23.2.4.1]{AndreIntAuxMotive}, for abelian varieties defined over $\Qbar$, there is an equality between $\trdeg \Q (\bt)$ (for the field generated by the coefficients of the matrix $\bt$) and the smallest dimension of Shimura subvariety containing it. More precisely, in \cite{FonsHigherRamanRelations}, special subvarieties also appear as obstructions to algebraic independence results (and conjecturally, the only ones), see \cite[Conjecture 13.1.2]{FonsHigherRamanRelations} for Hilbert modular surfaces and Hirzebruch-Zagier divisors. 

In very broad terms, in the higher dimensional setting the CM points are understood as special points, or zero dimensional special subvarieties. Because $\A_2$ (and any $\A_g$ for $g \geq 2$) does have positive dimensional special subvarieties, one would be inclined to think that special subvarieties should play a role in this generalization.

The answer we give to Question 1 is positive, and for Question 2, we can associate to every CM point a special subvariety of $\A_2$ where this type of exception for the transcendence degree \\$\trqjq$ is achieved for other points, not necessarily CM, in a "generic" sense, cf. Theorem \ref{TeoBigTeo}. This is a marked contrast with the statement of Theorem \ref{TeoCohenShigaWolfart}, which implies that for problems in terms of $\trdeg \Q(\t_1, \t_2, \t_3, j_1(\bt), j_2(\bt), j_3(\bt))$, the CM points are "isolated" exceptions. The computation of the transcendence degree $\trqjq$ for some CM points is conditional to Schanuel's conjecture \ref{ConjSchanuels}, or more precisely to the Gelfond-Schneider conjecture \ref{ConjGelfSchne}. 

It is noteworthy to highlight the role played by the Humbert singular relations (see Section \ref{sec:SubsectioHSR}), and that,  \emph{unconditionally}, knowledge about Humbert singular relations in moduli spaces of abelian surfaces gives insight on these questions stemming from the transcendence framework. We believe, however, that it might be difficult to generalize these results to higher genus curves: the antisymmetric relations in Remark \ref{RemMatricialHumbertAntisym} in the case of \emph{surfaces} lead to the elegant (and one dimensional) Humbert relations \eqref{EqHSR}, but this changes as soon as the abelian varieties have larger dimension.

\begin{rem} We state the result in terms of $\min_{\Sp_4(\Z)}$ because, as we will see in the following section (Example \ref{ExCurveCM5Root}), $\trqjq$ is not invariant under isomorphism of principally polarized abelian surfaces. With this notation, we mean 
\[
\min_{M \in \Sp_4(\Z)} \{\trdeg \Q (q_1(M\bt), q_2(M\bt), q_3(M\bt), j_{1}(M\bt)), j_2(M\bt)), j_3(M\bt)) \}.
\]
Remark that it is only the exponential functions $q_l = e^{2\pi i \tau_l}$ for $l=1,2,3$ that are not invariant under $\Sp_4(\Z)$.
\end{rem}

Throughout this paper, for $\tau \in \H_2$, by abuse of notation, we write $A_{\bt} \in \A_2$ instead of $[A_{\bt}] \in \A_2$ for the corresponding isomorphism class of ppas (principally polarized abelian surface). 

\begin{teo} Let $\bt \in \H_2$ a CM point with associated CM ppas $A_{\bt} \in \A_2^{ind}$, and set $q_l = e^{2\pi i \tau_l}$, $q = (q_1, q_2, q_3)$ and $j_l(q) = j_l(\bt)$, for $l=1,2,3$.
\begin{enumerate}
\item If $A_{\bt}$ is simple, it belongs to a Humbert surface of minimal discriminant, where, for any ppas defined over $\Qbar$ (which we associate with the notation $\tilde{q}$), it holds \\$\min_{\Sp_4(\Z)}\trqjqtil \leq 2$. Moreover, under Schanuel's conjecture \ref{ConjSchanuels}, $\min_{\Sp_4(\Z)}\trqjq = 2$.
\item If $A_{\bt} \simeq E' \times E''$ with $E',E''$ not isogenous CM elliptic curves, then the same statement holds verbatim. In this case, the special subvariety can also be taken as a special curve of type $\Q \times CM$. Moreover, under Schanuel's conjecture \ref{ConjSchanuels},\\ $\min_{\Sp_4(\Z)}\trqjq = 2$.
\item If $A_{\bt} \simeq E^2$ with $E$ a CM elliptic curve. Then $A_{\bt}$ belongs to one of the modular curves in the collection specified in Section \ref{sec:SecFamily}, and for any ppas in $C$ defined over $\Qbar$, it holds $\min_{\Sp_4(\Z)}\trqjqtil \leq 1$. \emph{Unconditionally}, \\$\trqjq = 1$.
\end{enumerate}
\label{TeoBigTeo}
\end{teo}

\subsection{Special subvarieties of \texorpdfstring{$\A_2$}{A2}}

The definition of special subvarieties requires the set-up of Shimura varieties. In $\A_2$, they may be explicitly described as loci of principally polarized abelian surfaces with some conditions (in particular, they are all of PEL type). We simply list them, as they appear in \cite[Table 1]{Daw-OrrUI}. We will eventually focus on the special subvarieties completely contained in $\A_2^{ind}$. We do not impose that a special subvariety needs to be connected.

The special surfaces are of the following two types, where the second one can be considered as a "degenerate case" of the first. They are jointly called Humbert surfaces. 
\begin{itemize}
\item Surfaces that parametrize ppas with real multiplication by a quadratic field.
\item Surfaces that parametrize ppas isogenous to a product of elliptic curves.
\end{itemize}

We will provide more details on these surfaces in Section \ref{sec:SubsectioHSR}. As an example, note that $\A_1 \times \A_1$ is a special surface. 

The special curves come in three types that we denote as follows.
\begin{itemize}
\item Type $\Q \times CM$: they parametrize ppas isogenous to a product $E \times E'$ with $E$ CM.
\item Shimura curves: they parametrize ppas with quaternionic multiplication (QM) by an indefinite division quaternion algebra over $\Q$.
\item Modular curves: they parametrize ppas isogenous to the square of an elliptic curve.
\end{itemize}

Modular curves can also be seen as "degeneration" of Shimura curves, with the quaternion algebra over $\Q$ taken as $\M_2(\Q).$ 

One can see, as an application of Theorem \ref{TeoBigThingHSR}, that points belonging to modular and Shimura curves always lie at the intersection of special surfaces, whereas a curve of type $\Q \times CM$ cannot be realized as an intersection of two special surfaces. Conversely, the intersection of two Humbert surfaces (when it is one dimensional) consists of a (finite union of) Shimura or modular curves. This perspective has been used several times in the literature to study said curves, for example, see \cite{HashiMuraba}, \cite{RunEndRingAbSur}, \cite{KaniGenHumbert} and \cite{BabGran}.

Moreover, for each of the special subvarieties, not all the three types of special points can occur:
\begin{itemize}
\item For the special surfaces parametrizing real multiplication, a CM point cannot be non-simple (by \cite[Corollary 2.7]{Goren}), so we find CM points which are simple or isotypic.
\item For the degenerate special surfaces, a CM point cannot be simple, by construction.
\item For the special curve of type $\Q \times CM$, a CM point cannot be simple, by construction.
\item For both modular and Shimura curves, a CM point can only be isotypic, via examination of the possible embeddings $B \hookrightarrow \End_0(A),$ for $B$ the indefinite quaternion algebra over $\Q$.
\end{itemize}

\begin{table}[h!]
\begin{tabular}{||c| c| c| c||} 
 \hline
 CM abelian surface & simple & non-simple & isotypic \\ [0.5ex] 
 \hline\hline
 $\Hum_{\Delta}$ & & NO  &  \\ 
 \hline
 $\Hum_{\delta^2}$ & NO &  &  \\
 \hline
 curve of type $\Q \times CM$ & NO & &  \\
 \hline
 Shimura curve & NO & NO &  \\
 \hline
 Modular curve & NO & NO &  \\ [1ex] 
 \hline
\end{tabular}
\caption{CM points on special subvarieties.}
\end{table}
 
\newpage
\subsection{Organization} 
The organization of this paper is as follows. In Section 2 we reduce our main problem to study linear relations between the coefficients of a period matrix $\bt \in \H_2$ of a ppas, and for CM points that is sufficient under the Gelfond-Schneider conjecture \ref{ConjGelfSchne}. In Section 3 we exemplify that Theorem \ref{TeoBigTeo} needs to be stated in terms of the $\Sp_4(\Z)$-orbit, because linear relations between the coefficients of $\bt$ are not invariant under the $\Sp_4(\Z)$-action, and give the minimum value of $\trqjq$ for CM points. In Section 4 we present Humbert singular relations, or linear relations involving the coefficients of $\bt$ \emph{and} $\det\bt$, which behave better for this action. We introduce the lattice of singular relations and the lattice of singular linear relations $\L_{\bt}$ and $\Llin_{\bt}$, which are positive definite lattices. We study $\L_{\bt}$ and, in particular, compute $\rank \L_{\bt}$, which is determined by $\End_0(A_{\bt})$. Section 5 is devoted to study $\rank \Llin_{\bt}$ and prove the main result Theorem \ref{TeoBigTeo}. The first two items of Theorem \ref{TeoBigTeo} are proved first, thanks to an application of Humbert's lemma Proposition \ref{PropHumbertsLemma}. The last item requires distinguishing between Shimura and modular curves: both have $\rank \L = 2$ generically. We study modular embeddings for both curves, searching for $\rank \Llin_{\bt} = 2$ generically in the image. Modular curves in the collection described in Section \ref{sec:SecFamily}, the ones studied in \cite{KaniThemodulispace}, admit embeddings with such property (and are the ones we use to finish the proof of Theorem \ref{TeoBigTeo}), while the classical quaternionic embeddings in \cite{Hashimoto} for Shimura curves do not. We further conjecture that $\rank \Llin_{\bt} \not= 2$ for whole $\Sp_4(\Z)$-orbit of any given point of these curves. In the Appendix, we give a partial answer to that, via Hirzebruch-Zagier divisors in Hilbert modular surfaces.
\\

\textbf{Acknowledgements:} We thank Daniel Bertrand for numerous fruitful discussions during and after a stay of the author at IMJ-PRG during Spring 2024, and valuable comments on previous drafts.
We also thank Damien Robert for the conversations during a visit to Universit\'e de Bordeaux, Yunqing Tang for pointing us to Hirzebruch-Zagier divisors and Aurel Page for answering our questions on quaternion algebras.

We thank Fabien Pazuki for guidance and useful remarks throughout this whole project, and also Tim With Berland and Fadi Mezher for several office conversations. 
\section{On linear dependence relations}

Let us first answer the analogous question for the $j$-invariant. If $\tau \in \H$ is quadratic imaginary, then $j(\tau)$ is algebraic and $\trdeg \Q(q, j(q)) = \trdeg \Q(q)$. But by the Gelfond-Schneider theorem \cite[Chapter 3, Section 2, Theorem 3.1]{FeldNest}, it holds  $q = e^{2\pi i \tau} = (-1)^{2\tau} \not\in \Qbar$, as $2\tau \in \Qbar \setminus \Q.$ 

Suppose now that $\bt \in \H_2$ is a CM point. By Theorem \ref{TeoCohenShigaWolfart}, 
\[
\{\bt \in \H_2|\: \bt \text{ is a CM point}\} = \{ \bt \in \H_2 \cap \M_2(\Qbar), \; j_1(\bt), j_2(\bt), j_3(\bt) \in \Qbar \}.
\]
Therefore, $$\trdeg \Q(q_1, q_2, q_3, j_1(q), j_2(q), j_3(q)) = \trdeg \Q(q_1, q_2, q_3) = \trdeg \Q(e^{2\pi i \tau_1}, e^{2\pi i \tau_2}, e^{2\pi i \tau_3}).$$ We now state Schanuel's conjecture.

\begin{conj}\emph{\textbf{(Schanuel's conjecture \cite[Chapter 6, page 260]{FeldNest})  }}Suppose \\$x_1, \ldots , x_n \in \C$ are $\Q$-linearly independent, then $\trdeg \Q(x_1, \ldots , x_n, e^{x_1}, \ldots, e^{x_n}) \geq n$. \label{ConjSchanuels}
\end{conj}

Conditionally to this conjecture, we have the following result, which can be seen as a Lindemann-Weierstrass type of statement for $e^{i\pi (\cdot)}$. Likewise, it is also a particular case of the Gelfond-Schneider conjecture (for $\al = -1 = e^{i\pi}$), independently attributed to Schneider \cite{SchneiderConjecture} and Gelfond \cite[Chapter 6, Conjecture 1 page 259]{FeldNest}.
\begin{conj}[\textbf{The Gelfond-Schneider conjecture}] Consider $\al \in \Qbar$, $\al \not=0,1$ and $x_1, \dots ,x_n \in \Qbar$ such that $1, x_1, \dots ,x_n$ are $\Q$-linearly independent. Then $\al^{x_1}, \ldots ,\al^{x_n}$ are algebraically independent. \label{ConjGelfSchne}
\end{conj}

We have included the classical deduction from Schanuel's conjecture (in the particular case $\al = e^{i\pi}$) for the reader's convenience, although historically it is an older conjecture.

\begin{lema}[Under Conjecture \ref{ConjGelfSchne}] Assume $1, x_1, \dots x_n \in \Qbar$ are $\Q$-linearly independent. Then $e^{i2\pi x_1}, \ldots ,e^{i2\pi x_n}$ are algebraically independent.

More generally, let $r =\dim \spn_{\Q}(1, x_1, \dots x_n)$, then $\trdeg \Q (e^{i2\pi x_1}, \ldots ,e^{i2\pi x_n}) = r-1.$ \label{LemaLWforPiExp}
\end{lema}

\begin{rem} The $\Q$-linearly independence with $1$ cannot be omitted in the statement. If $x_1, \ldots ,x_n$ satisfy an equation $\sum_{i=1}^n a_i x_i = b$ with $a_i, b \in \Z$ then $q_l = e^{i2\pi x_l}$ for $l=1, \ldots , n$ solve a multiplicative dependence relation $q_1^{a_1} \cdots q_n^{a_n} = 1.$  \label{RemInHomLinDep}
\end{rem}

\begin{proof}
It follows that $i2\pi , i2\pi  x_1, \ldots, i2\pi x_n$ are $\Q$-linearly independent, so by Schanuel's conjecture,
\begin{align*}
n+1 & \leq \trdeg \Q(i2\pi , i2\pi  x_1, \ldots, i2\pi x_n, e^{i2\pi}, e^{i2\pi x_1}, \ldots , e^{i2\pi x_n}) =^{*}\\
&= \trdeg \Q(i\pi , \underbrace{x_1, \ldots, x_n, 1}_{\in \Qbar}, e^{i2\pi x_1}, \ldots , e^{i2\pi x_n}) = \trdeg \Q(i\pi, e^{i2\pi x_1}, \ldots , e^{i2\pi x_n})
\end{align*}
where the equality in (*) comes from $\Qbar(i2\pi , i2\pi  x_1, \ldots, i2\pi x_n) = \Qbar(i\pi, x_1, \dots x_n)$. Therefore, $i\pi, e^{i2\pi x_1}, \ldots , e^{i2\pi x_n}$ are algebraically independent, which in particular implies our statement.

The second part has already been proven for $r=n+1$. For general $r$, take $r$ linearly independent elements among $1, \ldots,  x_n$. Assume first that we can take $1=x_1$, then the claim follows as for the case $r = n+1$. Otherwise, assume $x_1, \dots , x_r$ linearly independent, but with $\sum_{i=1}^r a_i x_i = b$ for not all zero $a_1, \dots, a_r, b \in \Z$. Applying Schanuel's conjecture as above for $i2\pi x_1, \ldots , i2\pi x_r$ results in $\trdeg \Q (e^{i2\pi x_1}, \ldots , e^{i2 \pi x_r}) \geq r-1$, but as in Remark \ref{RemInHomLinDep}, the transcendence degree cannot be maximal, as there is a multiplicative dependence among $e^{i2\pi x_1}, \ldots , e^{i2 \pi x_r}$. Therefore, $\trdeg \Q (e^{i2\pi x_1}, \ldots , e^{i2 \pi x_r}) = r-1.$
\end{proof}

\begin{cor}[Under Conjecture \ref{ConjGelfSchne}] Let $\bt \in \H_2$ a CM point. Then 
\[
\trdeg \Q(q_1, q_2, q_3, j_1(q), j_2(q), j_3(q)) = \dim \spn_\Q (1, \tau_1, \tau_2, \tau_3) - 1. 
\] \label{CorConjeForCM} \end{cor}
\begin{rem} In the proof of Lemma \ref{LemaLWforPiExp}, it was explicitly used that $x_i \in \Qbar$. In the general case, $\dim \spn_\Q (1, \tau_1, \tau_2, \tau_3)$ will say nothing about $\trqjq$. But if we additionally assume that $j_l(q) \in \Qbar$, them $\dim \spn_\Q (1, \tau_1, \tau_2, \tau_3)$ gives an upper bound.
\end{rem}

Therefore, Question 1 has naturally led us to study non-homogeneous linear dependence relations between the coefficients of $\bt$, and for CM points this analysis is sufficient (conditionally to the Gelfond-Schneider conjecture, or Schanuel's conjecture). In any case, these non-homogeneous linear dependence relations produce multiplicative dependence relations between the $q_l$'s, and give some partial information about the transcendence degree we are interested in.

Our strategy for Theorem \ref{TeoBigTeo} is to prove that for every CM point, there exists a special subvariety of $\A_2$ containing it and for which $\dim \spn_{\Q}(1, \tau_1, \tau_2, \tau_3)$ is generically constant through it, and equal to the value attained at the CM point. Therefore, we attach to every CM point a subvariety such that the points in this subvariety parametrizing abelian surfaces defined over $\Qbar$ are at least "of the same type" of exception as the CM point.

\section{First obstructions and the simplest case}
We first show via an example that the answer to Question 1, as stated, unfortunately must be answered with "it depends". 

\begin{ex} One of the most "famous" genus 2 curves with CM is
\[
C:y^2 = x^6 -x.
\]
In \cite[page 92]{Bost-Mestre-Moret} a period matrix is computed. For $\xi = e^{2\pi i /5}$, it is given by
\[
 \begin{pmatrix}
\tau_1 & \tau_2 \\ 
\tau_2 & \tau_3 \end{pmatrix} :=
\begin{pmatrix}
-\xi^4 & \xi^2 + 1 \\
\xi^2 + 1 & \xi^2-\xi^3.
\end{pmatrix}
\] \label{ExCurveCM5Root}
\end{ex}
It is easy to check that the $\tau_l$'s are $\Q$-linearly independent, using that $-\xi^4 = 1 + \xi + \xi^2 + \xi^3$, and that $1, \xi, \xi^2, \xi^3$ are $\Q$-linearly independent. Therefore, by Corollary \ref{CorConjeForCM}, \\$\trdeg\Q(q_1, q_2, q_3, j_1(q), j_2(q), j_3(q)) = 3$, and hence this is not an exception to \eqref{EqGenBehav}.

On the other hand, an application of Humbert's lemma (Proposition \ref{PropHumbertsLemma}, or more precisely of the algorithm described in \cite[Proposition 4.5]{Birk-Wilh}), gives us the matrix
\[
M = \begin{pmatrix}
-1 & 0 & 0 & 0 \\
0 & 1 & 0 & 1 \\
-3 & 0 & -1 & 0\\
0 & 2 & 0 & 3 \\
\end{pmatrix} \in \Sp_4(\Z).
\]
Then $\bt' = M\bt$ belongs to the $\Sp_4(\Z)$-orbit of $\tau,$ so the associated ppas are isomorphic. One can compute that
\[
\bt' = M\bt = \begin{pmatrix}
\frac{2}{55}\xi + \frac{9}{55}\xi + \frac{4}{55}\xi^2 + \frac{1}{11}\xi^3 & \frac{1}{11}\xi + \frac{1}{55}\xi^2 + \frac{4}{55}\xi^3 \\
\frac{1}{11}\xi + \frac{1}{55}\xi^2 + \frac{4}{55}\xi^3 & \frac{2}{5} + \frac{4}{55}\xi^2 + \frac{3}{55}\xi^2 + \frac{1}{55}\xi^3
\end{pmatrix},
\]
which solves the linear equation
\[
-\tau'_1 + \tau_2' + \tau'_3 = 0.
\]

Here is the verification with Sage:
\newpage
\begin{verbatim}
k = CyclotomicField(5)
xi = k.gen()
tau = Matrix([[-xi^4, xi^2 + 1],[xi^2 + 1,xi^2-xi^3]])
A = Matrix([[-1,0],[0,1]])
B = Matrix([[0,0],[0,1]])
C = Matrix([[-3,0],[0,2]])
D = Matrix([[-1,0],[0,3]])
Om = (A*tau + B)*((C*tau + D).inverse())
-Om[0][0] + Om[0][1] + Om[1][1]
\end{verbatim}

Hence, by Corollary \ref{CorConjeForCM}, for $\bt'$ we do get an exception to \eqref{EqGenBehav}. Question 1 should then be reformulated to a statement for the whole orbit under the action of $\Sp_4(\Z)$ by fractional linear transformations.

On a different note, we will first consider the "most exceptional" type of exceptions we have found among the CM points (as in they minimize $\dim \spn_{\Q}(1, \tau_1, \tau_2, \tau_3)$, see Corollary \ref{CorConjeForCM}). Remark that by definition of $\H_2$, it follows that $\tau_1, \tau_3 \not\in \Q$, therefore 
\[
\dim \spn_{\Q}(1, \tau_1, \tau_2, \tau_3) \geq \dim \spn_{\Q}(1, \tau_1) = 2.
\]

\begin{defi}For a complex abelian variety $A \cong \faktor{\C^g}{\Lambda}$ of dimension $g$, by \emph{big period matrix} we mean a matrix $\Pi \in \M_{g \times 2g}(\C)$ constructed from a choice of basis vectors $\la_{1}, \ldots, \la_{2g} \in \Lambda$ in terms of another choice of basis vectors $e_1, \ldots, e_g \in \C^g$. If we set $\Pi = ( \Omega_1, \Omega_2 )$, then by a \emph{(small) period matrix} we mean $\Omega_1 \in \M_{g}(\C)$, after a choice of basis such that $\Omega_2 = I_g$. Given a big period matrix for $A$, we will say that $\Omega_2^{-1}\Omega_1$ is a period matrix. 
\end{defi}

Alternatively, for ppav (by \cite[paragraph before Proposition 8.1.1]{Birk-LangeCAV})) a period matrix comes from the choice of $\la_1, \ldots, \la_{2g} \in \Lambda$ a \emph{symplectic} basis with respect to the principal polarization in $A$, and setting the first $g$ vectors $\la_1,\ldots \la_g$ as a basis of $\C^g.$ 

\begin{lema} An abelian surface $A$ admits a big period matrix belonging to $\mathcal{M}_{2 \times 4}(K)$ for $K$ quadratic imaginary field if and only if $A$ is isogenous to $E^2$, with $E$ an elliptic curve with complex multiplication. \label{LemaSimplestPM} \end{lema}

We now prove Lemma \ref{LemaSimplestPM}, which we found as part of a longer exercise in \cite[Exercise 10, section 5.6, page 142]{Birk-LangeCAV}. It is stated for abelian varieties of any dimension, and it gives a larger chain of equivalences:
\begin{enumerate}
\item The Picard number of $A$ (\textit{i.e.} the rank of its N\'eron-Severi group $\NS(A)$) is maximal, $\rank \NS(A) = g^2$.
\item $A$ is isogenous to a $E^g$, $E$ elliptic curve with complex multiplication.
\item $A$ admits a big period matrix in $\M_{g \times 2g}(K)$, for $K$ quadratic imaginary field.
\item $A$ is isomorphic (as unpolarized abelian varieties) to a product of $E_1 \times \ldots \times E_g$, $E_i$ pairwise isogenous elliptic curves with complex multiplication
\end{enumerate}

\begin{proof} We are inspired by the proof of \cite[Corollary 10.6.3]{Birk-LangeCAV}. A big period matrix for $A \simeq E^2$ necessarily has the shape
\begin{equation}
\begin{pmatrix}
\tau & 1 & 0 & 0 \\
0 & 0 & \tau & 1
\end{pmatrix} R, \label{EqBigPeriodMatrix}
\end{equation}
for $R \in \M_4(\Q)$ and $\tau$ such that $E = \C / (\tau \: 1) \Z^2$. Hence
$\begin{pmatrix}
\tau & 1 & 0 & 0 \\
0 & 0 & \tau & 1
\end{pmatrix} \in \M_{2 \times 4}(\Q(\tau))$, with $\Q(\tau)$ quadratic imaginary, so the big period matrix \eqref{EqBigPeriodMatrix} of $A$ also belongs to $\Q(\tau)$.

On the other direction, if $A$ admits a big period matrix in $\M_{2 \times 4}(\Q(\al))$ with $\al$ quadratic imaginary, then can solve the matrix $R$ as
\[
\begin{pmatrix}
  r_0 \al + r_4 &  r_1 \al + r_5 &  r_2 \al + r_6 &  r_3 \al + r_7 \\
 r_8\al + r_{12} & r_9 \al + r_{13}& r_{10} \al + r_{14} &r_{11} \al + r_{15} \end{pmatrix} = \begin{pmatrix}
\al & 1 & 0 & 0 \\
0 & 0 & \al & 1
\end{pmatrix} \left(\begin{array}{rrrr}
r_{0} & r_{1} & r_{2} & r_{3} \\
r_{4} & r_{5} & r_{6} & r_{7} \\
r_{8} & r_{9} & r_{10} & r_{11} \\
r_{12} & r_{13} & r_{14} & r_{15}
\end{array}\right)
\]
We need to verify that $R$ defines an isogeny between $A$ and $E^2$, with $E = \C/ (\al \: 1)\Z^2$ (and it is enough to show it as complex tori). First, there exists $m \in \Z$ such that $mR \in \M_4(\Z)$. If we set $\Pi'$ for the big period matrix of $A$, and $\Pi := \begin{pmatrix}
\al & 1 & 0 & 0 \\
0 & 0 & \al & 1
\end{pmatrix}$, then it holds
\[
(m I_2) \Pi = \Pi' (mR),
\]
and by \cite[Equation (1.1) before Proposition 1.2.3]{Birk-LangeCAV}, this is the compatibility condition sufficient to have a homomorphism of complex tori $E^2 \to A$. It is an isogeny because it is surjective (its analytic representation is multiplication by $m$) and between tori of the same dimension.
\end{proof}

\begin{cor} Let $\bt \in \H_2$ such that if $A_{\bt} \simeq E^2$, with $E$ a CM elliptic curve. Then \\$\trqjq = 1.$ \label{CorCompEisotypic}
\end{cor}

\begin{proof}
This follows immediately from Lemma \ref{LemaSimplestPM}. If $A_{\bt} \simeq E^2$ as in the statement, then it admits a big period matrix $(\Omega \; \Omega')$ lying in $\mathcal{M}_{2 \times 4}(K)$ for $K$ a quadratic imaginary field. Therefore, $\bt \in \M_2(K)$, so $\dim \spn_\Q (1, \tau_1, \tau_2, \tau_3) \leq \dim_{\Q} K = 2$. Hence \\$\trqjq \leq 1$, and equality holds by the Gelfond-Schneider theorem \cite[Chapter 3, Section 2, Theorem 3.1]{FeldNest}, as $q_1 = e^{2\pi i \t_1}$ is transcendental.

We notice that $\bt' \in \M_2(K)$ along the whole $\Sp_4(\Z)$-orbit, because for $\begin{pmatrix}
A & B \\ C & D
\end{pmatrix} \in \Sp_4(\Z)$, it follows that $A\bt + B \in \M_2(K)$ and $(C\bt + D)^{-1} \in \M_2(K)$, as $K$ is a number field.
\end{proof}

\begin{rem} The proof of Corollary \ref{CorCompEisotypic} is thus unconditional, in particular is not relying on Schanuel's conjecture.
\end{rem}

\section{Humbert singular relations}
To go beyond in our study of linear relations between $\t_1, \t_2, \t_3$ we will know focus on studying linear relations between $1, \tau_1, \tau_2, \tau_3, \tau_2^2 - \tau_1 \tau_3$, \textit{i.e.} equations (over $\Z$) of the form
\begin{equation}
a\tau_1 + b\tau_2 + c\tau_3 + d(\tau_2^2 - \tau_1\tau_3) + e= 0, \label{EqHSR}
\end{equation}
which behave better under the action of $\Sp_4(\Z)$ (see Lemma \ref{LemmaRankLbt}). Returning to Example \ref{ExCurveCM5Root}, the coefficients of both $\bt$ and $\bt'$ satisfy such a relation, and the one for $\bt'$ happens to have $d=0.$
For $\bt$, as $\bt \in \M_2(\Q(\xi))$ with $\xi$ a primitive $5$-th root of unity, then $\dim \spn_{\Q}(1, \tau_1, \tau_2, \tau_3, \tau_2^2 - \tau_1 \tau_3) \leq  \dim \Q(\xi) = 4$, so there must exist a relation as in \eqref{EqHSR}.

\subsection{Basics on Humbert singular relations} \label{sec:SubsectioHSR}

The source material for this is \cite[section 4]{Birk-Wilh}, a modern exposition of \cite{Humbert}. Additionally, see \cite[Section 3]{Lin-Yang}. We will include the proofs of the relevant results for readability. 

For a ppas $A_{\bt}$, we consider two representations for $\End(A_{\bt})$. By \cite[Prop 1.2.1]{Birk-LangeCAV}, the endomorphism of the complex tori $A_{\bt} = \faktor{\C^2}{\Lambda_{\bt}}$ with $\Lambda_{\bt} = (\bt \: I_2) \Z^4$ correspond to $\C$-linear maps $F:\C^2 \to \C^2$ such that $F(\Lambda) \subset \Lambda$. Therefore, considering the $\C$-linear map $\C^2 \to \C^2$ induces the \emph{analytical representation} $\rho_{a,\bt}: \End(A_{\bt}) \to \mathcal{M}_2(\C)$, and considering the $\Z$-linear map on the lattice $\Lambda \to \Lambda$ induces the \emph{rational representation} $\rho_{r,\bt}: \End(A_{\bt}) \to \mathcal{M}_4(\Z)$. There is an extra compatibility condition: for any $f \in \End(A_{\bt})$, 
\[
\rho_{a,\bt}(f)(\Lambda) \subset \Lambda,
\]
or, in terms of matrices, for any $f \in \End(A_{\bt})$,
\begin{equation}
\rho_{a,\bt}(f) (\bt \: I_2) = (\bt \: I_2) \rho_{r,\bt}(f). \label{EqRelAnRat}
\end{equation}
Conversely, any pair of matrices in $\M_4(\Z)$ and $\M_2(\C)$ satisfying \eqref{EqRelAnRat} induces an element in $\End(A_{\bt})$.

\begin{rem} Notice that $\rho_{r,\bt}$ is not surjective. Setting $\rho_{r,\bt}(f) = \begin{pmatrix} A & B \\ C & D \end{pmatrix}$ for matrices in $\mathcal{M}_2(\Z)$, this is equivalent to
\[
\left( \rho_{a, \bt}(f) \bt \quad \rho_{a, \bt}(f) \right) = \left( \bt A + C \quad \bt B  + D \right),
\]
and solving for $\rho_{a, \bt}(f) \bt$, leads to 
\begin{equation}
\left( \bt A + C \right) = \left( \bt B + D \right) \bt. \label{EqRelForEnd}
\end{equation}

Remark that as $A_{\bt}$ always admits a group structure, $\Z \subset \End(A_{\bt})$, and via the rational representation $\Z$ is identified with matrices $n I_4$ with $n \in \Z$. Note that for $f \in \Z$, \eqref{EqRelForEnd} is trivially satisfied as $A=D=nI_2$ and $B=C=0$.

This also illustrates the fact that, generically in $\A_2$, we have $\End(A) = \Z$.
\end{rem}

These representations depend on $\bt \in \H_2$. For $M \in \Sp_4(\Z)$, then
\begin{equation}
\rho_{r, M\bt} ={}^{t}\!M^{-1} \rho_{r,\bt}{}^{t}\!M, \label{EqRationalRepUnderSp4}
\end{equation}
as ${}^{t}\!M$ is the rational representation on the isomorphism from $A_{\bt}$ to $A_{\bt'}$. For $\rho_{a,\bt}$, if we set $M = \begin{pmatrix} A & B \\ C & D \end{pmatrix}$, the analytic representation of said isomorphism is ${}^{t}\!(C \tau + D)$, hence
\begin{equation}
\rho_{a, M\bt} = {}^{t}\!(C \tau + D)^{-1} \rho_{a,\bt} {}^{t}\!(C \tau + D). \label{EqAnalRepUnderSp4}
\end{equation}

The principal polarization induces an anti-involution on $\End(A)$ called the \emph{Rosati involution}, which we denote $':\End(A) \to \End(A)$, corresponding to the adjoint operator to the Riemann form $E = \begin{pmatrix}
0 & I \\ -I & 0
\end{pmatrix}$. That implies (\cite[Proposition 5.1.1]{Birk-LangeCAV}) that for all $\lambda, \mu \in \Lambda$ and $f \in \End(A_{\bt}),$
\[
E(\rho_{r,\bt}(f)(\lambda), \mu) = E(\lambda, \rho_{r,\bt}(f')(\mu)),
\]
and in terms of matrices,
\begin{equation}
{}^t\!\rho_{r,\bt}(f)\begin{pmatrix}
0 & I_2 \\
-I_2 & 0
\end{pmatrix} = \begin{pmatrix}
0 & I_2 \\
-I_2 & 0
\end{pmatrix}
\rho_{r,\bt}(f'). \label{EqSymRatRep}
\end{equation}
Likewise, for the analytic representation, as $H =(\Im \bt)^{-1}$, we have
\begin{equation}
\overline{\rho_{a, \bt}(f')} \Im \bt = \Im \bt \, {}^t\!\rho_{a, \bt}(f), \label{EqSymAnRep}
\end{equation}
where $\bar{\cdot}$ indicates complex conjugation.
We call $f \in \End(A)$ \emph{symmetric} if $f' = f$. We set $\End^s(A) \subset \End(A)$ the subring of symmetric endomorphisms. We have the following result (from \cite[Lemma 4.1]{Birk-Wilh}) for symmetric endomorphisms.

\begin{lema} Consider an endomorphism $f\in \End(A_{\bt})$, with rational representation given by $\rho_{r,\bt}(f) = \begin{pmatrix}
A & B \\C& D
\end{pmatrix}$. Denote $A = \begin{pmatrix}
a_1 & a_2 \\ a_3 & a_4
\end{pmatrix}$, $b = B_{12}$, and $c = 	C_{12}$. For the following statements, we have $(1) \iff (2)$, and any of them implies $(3)$.
\begin{enumerate}
\item The endomorphism $f$ is symmetric.
\item The matrices $B$ and $D$ are antisymmetric and $D = {}^t\!A$: more precisely it follows that
\[
B = \begin{pmatrix}
0 & b \\
-b & 0
\end{pmatrix}, \, C= \begin{pmatrix}
0 & c \\
-c & 0
\end{pmatrix}, \, D={}^t\!A.
\]

\item If $\bt = \begin{pmatrix}
\t_1 & \t_2 \\ \t_2 & \t_3
\end{pmatrix}$, then it satisfies
\[
a_2 \t_1 + (a_4 - a_1)\t_2 - a_3 \t_3 + b\left( \t_2^2 - \t_1\t_3\right) + c = 0.
\]
\end{enumerate}
In particular, for non-trivial equations, $\{ 1, \t_1, \t_2, \t_3, \det(\bt) \}$ are $\Q$-linearly dependent.
 \label{LemaHumertRelations}
\end{lema}
\begin{proof}
The equivalence $(1) \iff (2)$ comes from unwinding \eqref{EqSymRatRep}, which gives
\[
\begin{pmatrix}
-{}^t\!C & {}^t\!A \\-{}^t\!D & {}^t\!B 
\end{pmatrix} = \begin{pmatrix}
C & D \\ -A & -B
\end{pmatrix},
\]
hence $D = {}^t\!A$ and $B$ and $C$ are antisymmetric. This equation is equivalent to $f = f'.$

Now $(2) \implies (3)$ is a consequence of \eqref{EqRelForEnd}, that is equivalent to
\begin{align*}
\left( \bt A + C \right) &= \left( \bt B +  {}^t\!A\right) \bt, \\
\bt B \bt - C &= (\bt A - {}^t\!(\bt A)), 
\end{align*}
and as the three matrices involved are antisymmetric, the matrix equation above holds if and only if it holds for the off diagonal element, which is the equation of the statement.
\end{proof}

\begin{rem} As before, $\Z \subset \End(A)$ and actually $\Z \subset \End^s(A)$, but this yields trivial equations for $\bt.$ In addition, note that if $f \in \End^s(A)$, then $n + m f \in \End^s(A)$, with rational representation $nI_4 + m\rho_{r,\bt}(f).$ \label{RemHumbRelation}
\end{rem}

Remark that on the linear equation in Lemma \ref{LemaHumertRelations} we can take $a_1=0$, as that means considering $f- a_1 \in \End_s(A_{\bt})$ instead of $f$, with rational representation $\rho_{r,\bt}(f) - a_1 I_4$.

For given $(a,b,c,d,e)\in \Z^5$, one can define the following matrix $R_{(a,b,c,d,e)} = \begin{pmatrix}
A & B \\ C & {}^t\!A
\end{pmatrix}$
where
\[
A = \begin{pmatrix} 0 & a \\
-c & b
\end{pmatrix}, \; B = \begin{pmatrix}
0 & d \\
-d & 0
\end{pmatrix}, \, C= \begin{pmatrix}
0 & e \\
-e & 0
\end{pmatrix}.
\]

\begin{lema} If we have $(a,b,c,d,e) \in \Z^5$ such that
\[
a\t_1 + b\t_2 + c\t_3 + d(\t^2 - \t_1\t_3) + e= 0,
\]
then there exists $f \in \End^s(A_{\bt})$ with $\rho_{r,\bt}(f) = R_{(a,b,c,d,e)}$. \label{LemaHSRtoSymEnd}
\end{lema}
\begin{proof}
By the argument at the end of the proof of Lemma \ref{LemaHumertRelations}, the equation in the statement is in fact equivalent to
\[
\left( \bt A + C \right) = \left( \bt B +  {}^t\!A\right) \bt,
\]
which is equivalent to \eqref{EqRelForEnd}, and to the compatibility condition \eqref{EqRelAnRat}. Hence, there exists an endomorphism $f \in \End(A_{\bt})$ with such a rational representation. It is symmetric by Lemma \ref{LemaHumertRelations}.
\end{proof}

\begin{rem} In the proof of Lemma \ref{LemaHumertRelations} and Lemma \ref{LemaHSRtoSymEnd}, the only instance where we use that the abelian variety has dimension $g=2$ is that the matrix equation involves antisymmetric matrices in $\M_2(\C)$, so it is a subspace of dimension one. That is why Humbert singular relations only arise naturally in this dimension. \label{RemMatricialHumbertAntisym}
\end{rem}

In conclusion, the existence of $f \in \End^s(A_{\bt})$ with $f \not\in \Z$ is determined, via $\rho_{r,\bt}(f)$, by certain quadratic relations for the components of the period matrix $\bt$. 

\begin{defi} We call \emph{Humbert singular relation} (HSR for short) the following equation for $a,b,c,d,e \in \Z$,
\[
a\t_1 + b\t_2 + c\t_3 + d(\t_2^2 - \t_1\t_3) + e =0.
\]
If in addition $\gcd(a,b,c,d,e) = 1$, we say that the relation is \emph{primitive}.
\end{defi}
Let us consider now the analytic representation $\rho_{a,\bt}(f)$ of a symmetric endomorphism, and its characteristic polynomial as a matrix in $\M_2(\C)$:
\[
P_{f}(t):= \det(tI_2 - \rho_{a,\bt}(f)) = t^2 - \tr(\rho_{a,\bt}(f))t + \det(\rho_{a,\bt}(f)).
\]
By \eqref{EqRelAnRat} and Lemma \ref{LemaHumertRelations}, for $f \in \End^s(A_{\bt})$,
\[
\rho_{a,\bt}(f) = \bt B +{}^t\!A = \begin{pmatrix}
-d\tau_2 & d\tau_1 - c \\ -d\tau_3 + a & d\tau_2 + b
\end{pmatrix}.
\]
Hence, it follows (using the singular relation solved by $\bt$):
\begin{align*}
\tr(\rho_{a,\bt}(f)) & = b, \\
\det(\rho_{a,\bt}(f)) &= ac + de.
\end{align*}
Therefore, $P_f(t) = t^2 - bt + (ac + de)$. More is true, $\rho_{a,\bt}$ induces a ring isomorphism (\cite[Corollary 4.4]{Birk-Wilh})
\[
\{ n + mf, \: n,m\in \Z \} \to \faktor{\Z[t]}{\left(t^2 -bt + ac + de\right)}.
\]
Note that $\Delta = \Delta(f) := \Disc(P_f) = b^2 - 4(ac + de)$. Remark that $\Delta \equiv 0,1$ (mod $4$).

\begin{lema} If $f \in \End^s(A_{\bt})$ then $\Delta(f) \geq 0$, and $\Delta (f) = 0$ if and only if $f \in \Z$. \label{LemaDiscriSymEndom}
\end{lema}

\begin{proof}
(From \cite[Proposition 4.7]{Birk-Wilh}). By \eqref{EqSymAnRep}, it follows that $\rho_{a,\bt}(f) \Im \bt$ is a Hermitian matrix. Consider $\lambda_k$ an eigenvalue of $\rho_a(f)$ with eigenvector $v_k$, for $k=1,2$. Because $f$ is symmetric,
\[
\lambda_k H(v_k, v_k) = H(\rho_a(f) v_k, v_k) =H(v_k, \rho_a(f) v_k) = \overline{\lambda_k} H(v_k, v_k),
\]
therefore $\lambda_k$ is real, and is $\lambda_k$ is a root of $P_f$, it follows that $\Delta(f) \geq 0$ and $\Delta(f) = (\lambda_1 - \lambda_2)^2$. It can only be $0$ if the two eigenvalues of $\rho_a(f)$ coincide and are integers, say $m \in \Z$, as $P_f(t) \in \Z[t]$. Therefore, $\rho_a(f) - mI_2 = \rho_a(f - m)$ is a nilpotent matrix (with exponent two, necessarily, for it is a $2 \times 2$ matrix). One can check that the following linear algebra problem on $\M_2(\C)$
\[
\begin{cases}
Z^2 = 0\\
Z\Im(\bt) \text{ is Hermitian}
\end{cases}
\]
only has as solution $Z = 0$, and in particular for us $\rho_a(f - m) = 0$, and $f\in \Z$. A possible proof goes by inspection: in dimension two a nilpotent matrix is characterized by trace and determinant zero, so there exists $z,u,v\in \C$ with $-z^2 = uv$ such that
\[
Z = \begin{pmatrix}
z & u \\ v & -z
\end{pmatrix}.
\]
One can impose the Hermitian conditions on
\[
\begin{pmatrix}
z & u \\ v & -z
\end{pmatrix} \begin{pmatrix}
\Im (\t_1) & \Im (\t_2) \\ \Im (\t_2) & \Im (\t_3) 
\end{pmatrix} = \begin{pmatrix}
z\Im (\t_1) + u \Im (\t_2) & z\Im (\t_2) + u \Im (\t_3) \\ v\Im (\t_1) -z\Im (\t_2) & v\Im (\t_2) - z\Im (\t_3)
\end{pmatrix},
\]
resulting on the equations (where we set $y_i = \Im(\t_i)$, remark that $y_1, y_3 > 0$):
\begin{align*}
zy_1 + u y_2 &= \bar{z}y_1 + \bar{u}y_2,\\
vy_2 - zy_3 &= \bar{v}y_2 - \bar{z}y_3, \\
vy_1 - zy_2 &= \bar{z}y_1 + \bar{u}y_3,
\end{align*}
equivalent to
\begin{align}
\Im(z)y_1 &= -\Im(u)y_2, \label{Eq1} \\
\Im(z)y_3 &= \Im(v)y_2,  \label{Eq2} \\
v &= \frac{2\Re(z)y_2 + \bar{u}y_3 }{y_1}\label{Eq3}.
\end{align}
Note that \eqref{Eq3} implies 
\begin{equation}
y_1 \Im(v) = -y_3\Im(u). \label{Eq3'}
\end{equation}
As $-z^2 = uv$, it follows from \eqref{Eq3} that
\begin{equation}
y_1 z^2 + 2y_2u\Re(z) + y_3 |u|^2 = 0, \label{Eq4}
\end{equation}
which implies, as $\Re(z^2) = \Re(z)^2 - \Im(z)^2$
\[
y_1(\Re(z)^2 - \Im(z)^2) + 2y_2\Re(u)\Re(z) + y_3(\Re(u)^2 + \Im(u)^2) = 0,
\]
or equivalently,
\begin{equation}
y_1\Re(z)^2 + 2y_2\Re(u)\Re(z) +y_3\Re(u)^2 = y_1\Im(z)^2 - y_3\Im(u)^2. \label{Eq5}
\end{equation}
We rewrite the right hand side. Together with \eqref{Eq1} and \eqref{Eq2},
\begin{align*}
y_1\Im(z)^2 - y_3\Im(u)^2 &= -\Im(z)\Im(u)y_2 + y_1\Im(v) \Im(u) = \Im(u) \left( -y_2\Im(z) + y_1\Im(v) \right) \\
&= \Im(u) \left( \frac{-\Im(v)y_2^2}{y_3} + y_1\Im(v) \right) = \Im(u) \Im(v) \frac{y_1 y_3 - y_2^2}{y_3} \leq 0,
\end{align*}
because $\Im(u)\Im(v) \leq 0$ by \eqref{Eq3'}, $\det(\Im \bt) > 0$ and $y_3 > 0$. But on the other hand, 
\[
y_1\Re(z)^2 + 2y_2\Re(u)\Re(z) +y_3\Re(u)^2 = (\Re(z) \; \Re(u)) \begin{pmatrix}
y_1 & y_2 \\ y_2 & y_3
\end{pmatrix}
\begin{pmatrix}
\Re(z) \\ \Re(u)
\end{pmatrix} \geq 0,
\]
as $\Im(\bt)$ is positive definite. Therefore, the only possibility for Equation \eqref{Eq5} to hold is if $\Re(z) = \Re(u) = 0$ and $\Im(u) \Im(v) = 0$. By \eqref{Eq3} we have $0 =\Im(u) = \Im(v)$, and by \eqref{Eq1} $\Im(z) = 0$. Hence $z,u,v \in \R$ and we have $z = u = 0$, so by \eqref{Eq3} we conclude $v = 0$.
\end{proof}

\begin{rem} We remark that the end of the proof of Lemma \ref{LemaDiscriSymEndom} implies that there are no \emph{symmetric} nilpotent endomorphisms. When $A_{\bt} \simeq E^2$, the endomorphism algebra is a matrix algebra and contains nilpotent elements, therefore they cannot be symmetric with respect to the Rosati involution. 

Compare this with \cite[Section 5.3]{Birk-LangeCAV} (for any dimension $g$): there exist symmetric \emph{idempotents}, the norm-endomorphisms detecting abelian subvarieties. Roughly speaking, they are constructed as "projections" $f_X = f \in \End^s(A)$ to a subvariety $X \subset A$ such that \\$f \lvert_{X}$ $= e(X)_{[X]}$, for $e(X)$ the exponent of the induced polarization. These endomorphisms are symmetric with respect to the Rosati involution and satisfy $f^2 = e(X)f$. One can further consider the (canonical) idempotent $e(X)^{-1}f \in \End^s_0(A)$. By \cite[Theorem 5.3.2]{Birk-LangeCAV}, the symmetric idempotents of $\End_0(A)$ are in bijection with the abelian subvarieties of $A.$ 

In particular, for ppas, it follows that for $f \in \End^s(A)$, then either $f$ is invertible in $\End_0(A)$ or there exists $\mu \in \Z$ such that $f-\mu$ is a multiple of an idempotent element in $\End_0(A)$.
\end{rem}

After Lemma \ref{LemaDiscriSymEndom}, we can consider the following definition. Remark that $\Hum_0 = \A_2$, so we do not consider it.

\begin{defi} For $\Delta > 0 $, $\Delta \equiv 0,1$ (mod $4$), we call \emph{Humbert surface of discriminant} $\Delta$ the following locus in $\A_2$:
\[
\Hum_{\Delta} := \{A \in \A_2| \text{ there exists } f\in \End^s(A) \text{ with } \Delta(f)= \Delta \} \subset \A_2.
\]
\end{defi}

Note that $\Hum_{m^2 \Delta}$ is well defined for any $m\in \Z$ and $\Hum_{\Delta} \subset \Hum_{m^2 \Delta}.$ 

We will see in Lemma \ref{LemmaRankLbt} in the following section that the discriminant of a HSR is invariant under $\Sp_4(\Z)$-action, so its preimage in $\H_2$ is
\begin{gather*}
\{\bt \in \H_2: \text{ there exists a HSR } l \text{ with } \Delta(l) = \Delta \text{ such that } \bt \text{ satisfies } l \}.
\end{gather*}
It follows that imposing a \emph{primitive} HSR minimizes the discriminant.

We distinguish the following cases in terms of this discriminant, which is the content of \cite[Proposition 4.8, Proposition 4.9]{Birk-Wilh}.

\begin{prop}
Suppose that $\bt$ satisfies a primitive HSR with discriminant $\Delta$. Then the following holds:
\begin{itemize}
\item If $\Delta$ is a fundamental discriminant, then there exists an embedding of the ring of integers
\[
\O_{\Q(\sqrt{\Delta})} \hookrightarrow \End^s(A).
\]
More generally, if $\Delta$ is the discriminant of an order $\O$ in a real quadratic field, then there exists an embedding
\[
\O \hookrightarrow \End^s(A).
\]
In both cases, there is an embedding of a real quadratic field $K \hookrightarrow \End^s_0(A)$.

\item If $\Delta = \delta^2$, then $A$ is not simple. More precisely, there exists an isogeny of degree $\delta^2$:
\[
\left( E_1 \times E_2, p_1^* \O_{E_1}(\delta) \otimes p_2^* \O_{E_2}(\delta)  \right) \to (A, H).
\]
\end{itemize} \label{PropDeltasAndQuadFieldIsogenies}
\end{prop}

A diagonalization by blocks algorithm (see \cite[Proposition 4.5]{Birk-Wilh}) for $\rho_{r,\tau}(f)$ allows us to choose a basis in order to obtain a symmetric (non trivial) endomorphism of the following special shape:
\[
\begin{pmatrix}
A & 0 \\ 0 & {}^t\!A
\end{pmatrix}, \: A= \begin{pmatrix}
0 & a \\ -1 & b
\end{pmatrix}
\]
with 
\[
(a,b) = \begin{cases} \left(\frac{-\Delta}{4}, 0\right) & \text{ if } \Delta \equiv 0 \text{ mod } 4, \\
\left(\frac{1 - \Delta}{4}, 1 \right) & \text{ if } \Delta \equiv 1 \text{ mod } 4.
\end{cases}
\]
Therefore, as per Lemma \ref{LemaHumertRelations}, given $\bt \in \H_2$ satisfying a HSR of discriminant $\Delta$, we have a distinguished one, which we call the \emph{normalized} HSR of discriminant $\Delta.$

\begin{prop}[Humbert's lemma] For $\bt \in \H_2$, if $A_{\bt}$ satisfies a Humbert singular relation of discriminant $\Delta$, then there exists $M \in \Sp_4(\Z)$ such that $\bt' = M\bt$ solves:
\begin{align*}
\frac{-\Delta}{4}\tau_1' + \tau_3' = 0 & \text{ if } \Delta \equiv 0 \text{ mod } 4,\\
\frac{1 - \Delta}{4} \tau_1' + \tau_2' + \tau_3' = 0 & \text{ if } \Delta \equiv 1 \text{ mod } 4.
\end{align*} \label{PropHumbertsLemma}
\end{prop}
We will say more on how Humbert singular relations behave under $\Sp_4(\Z)$-action in the next section.

Finally, we finish this section with the relation between $\End^s(A)$ and the N\'eron-Severi group $\NS(A)$, the group of divisors modulo algebraic equivalence, or the image of the first Chern class $\Pic(A) \to H^2(A,\Z)$. For a line bundle $L$ on $A$ and for any $x \in A$, the line bundle $t_x^*L \otimes L^{-1}$ has first Chern class zero, so it belongs to $\Pic^0(A) = \hat{A}$. This induces a map $\phi_L: A \to \hat{A}$ given by $x \mapsto t_x^*L \otimes L^{-1}$, that only depends on the first Chern class of $L$ in $\NS(A)$. For a principal polarization $H_0$, it is moreover an isomorphism $\phi_{0}: A \to \hat{A}$. We can consider the assignment:
\begin{align*}
\NS(A) & \to \End(A) \\
L & \mapsto \phi_{0}^{-1} \circ \phi_{L}.
\end{align*}
In fact, it is compatible with the Rosati involution, restricting to $\NS(A) \to \End^s(A)$. More can be said:

\begin{prop} The assignment
\begin{align*}
\NS(A) & \to \End^s(A) \\
L & \mapsto \phi_{0}^{-1} \circ \phi_{L}
\end{align*}
is a group isomorphism. \label{PropIsomNSEnds}
\end{prop}
\begin{proof}
This is \cite[Proposition 5.2.1]{Birk-LangeCAV}.
\end{proof}
This actually generalizes to an isomorphism of $\Q$-vector spaces for $H_0$ not necessarily principal. Conversely, we can think of $\End^s(A)$ as parametrizing other line bundles in $A$, in particular, other (principal) polarizations: by \cite[Theorem 5.2.4]{Birk-LangeCAV}, there is a bijection between the polarizations of degree $d$ in $\NS(A)$ and the totally positive endomorphisms in $\End^s(A)$ of norm $d$.

\subsection{The lattice of singular relations}

\begin{defi} For $\bt \in \H_2$, we denote \emph{the lattice of Humbert singular relations} the following $\Z$-module  \label{DefLattoceSingularRelations}
\[
\L_{\bt} :=  \{ (a,b,c,d,e) \in \Z^{5}:  a\tau_1 + b\tau_2 + c\tau_3 + d(\tau_2^2 - \tau_1\tau_3) + e = 0\},
\]
equipped with the positive definite (by Lemma \ref{LemaDiscriSymEndom}) quadratic form induced by the discriminant $\Delta$. It forms a positive definite integral lattice (as in finite rank free abelian group with a symmetric bilinear form).

Likewise, we consider the sublattice $\Llin_{\bt}:= \L_{\bt} \cap \{ d=0\}$, or equivalently, the sublattice spanned by singular relations that do not involve $\tau^2 - \tau_1\tau_3$. Therefore, one can consider $\L_{\bt} \otimes \Q \subset \Q^5$ and $\Llin_{\bt} \otimes \Q \subset \Q^4.$
\end{defi}

We are ultimately interested in the rank of $\Llin_{\bt}$. We will study $\rank \L_{\bt}$ and then compare it with $\rank \Llin_{\bt}.$ 

Some of these results admit alternative proofs. First, by Lemma \ref{LemaHumertRelations}, Remark \ref{RemHumbRelation} and Lemma \ref{LemaHSRtoSymEnd}, there is a short exact sequence of abelian groups:
\begin{equation}
0 \to \Z \to \End^s(A_{\bt}) \to \L_{\bt} \to 0, \label{EqSESNeronSeveri}
\end{equation}
where the homomorphism $\End^s(A_{\bt}) \to \L_{\bt}$ is given by 
\[
f \mapsto \rho_{r,\bt}(f) = \begin{pmatrix}
\begin{pmatrix}
a_1 & a_2 \\ a_3 &a_4
\end{pmatrix} & \begin{pmatrix}
0 & b \\ -b & 0
\end{pmatrix} \\ \begin{pmatrix}
0 & c \\ -c & 0
\end{pmatrix} & {}^t\!\begin{pmatrix}
a_1 & a_2 \\ a_3 &a_4
\end{pmatrix}
\end{pmatrix} \mapsto (a_2, (a_4 - a_1), -a_3, b, c)
\]
If one combines \eqref{EqSESNeronSeveri} and Proposition \ref{PropIsomNSEnds}, it follows that $\rank \L_{\bt}= \rank{\NS(A_{\bt})}-1 $. Hence, one can use results for $\rank\NS(A_{\bt})$ instead. However, we believe that a proof via explicit manipulation of the lattices is insightful in our setting, for that we are also able to understand the effect on the quadratic form $\Delta$.

In that direction, there is a more intrinsic definition of $(\L_{\bt}, \Delta)$ in $\faktor{\NS(A_{\bt})}{\Z H_0}$, with the quadratic form defined in terms of the intersection pairing, see \cite[Section 2]{KaniThemodulispace}.

\begin{rem} The following is true:
\begin{itemize} \item We have the following short exact sequences of $\Q$-vector spaces
\[
0 \mapsto \L_{\bt} \otimes \Q \mapsto \Q^5 \mapsto \spbt \mapsto 0,
\]
and
\[
0 \mapsto \Llin_{\bt} \otimes \Q \mapsto \Q^4 \mapsto \spn_\Q(1,\t_1, \t_2, \t_3) \mapsto 0.
\]
Because $\tau_1 \not\in \Q$, the dimension of both spans is at least $2$. Consequently $0 \leq \rank  \L_{\bt} \leq 5-2 = 3$, and $0 \leq \rank \L_{\bt}^{lin} \leq 4-2=2.$
\end{itemize}
\end{rem}

The extremal bounds for $\rank \L_{\bt}$ are realized.

\begin{lema} Let $\bt \in  \H_2$. Then
\begin{itemize}
\item $\rank \L_{\bt} = 0$ if and only if $\End_0(A_{\bt}) = \Q$.
\item $\rank \L_{\bt} = 3$ if and only if $A_{\bt} \simeq E^2$ for $E$ a CM elliptic curve. \label{LemaRankLtExtremeCases}
\end{itemize} 
\end{lema}

\begin{proof}
The case $\rank \L_{\bt} = 0$ corresponds to the absence of any Humbert singular relation whatsoever, hence, $\End^s(A_{\bt}) = \Z$ by Lemma \ref{LemaHumertRelations} and Lemma \ref{LemaHSRtoSymEnd}. By \cite[Theorem 5.3.2]{Birk-LangeCAV}, any abelian subvariety is associated to a symmetric idempotent, so $A_{\bt}$ is necessarily simple if $\End_0^s(A_{\bt}) = \Q$. By Albert's classification of endomorphisms algebra of simple abelian varieties \cite[Proposition 5.5.7]{Birk-LangeCAV}, and the fact that the case III does not occur for abelian surfaces, it follows that $\End_0(A_{\bt}) = \Q.$

For the second part, note that $\rank \L_{\bt} = 3$ is equivalent to $\dim \spn_{\Q} (1, \tau_1, \tau_2, \tau_3, \tau_2^2 - \tau_1\tau_3)=2$ and hence,
\[
\spn_{\Q} (1, \tau_1) = \spn_{\Q} (1, \tau_1, \tau_2, \tau_3) = \spn_{\Q} (1, \tau_1, \tau_2, \tau_3, \tau_2^2 - \tau_1\tau_3). 
\]
Hence, there exists $a,b ,a', b', a'', b'' \in \Q$ such that:
\begin{align*}
\tau_2 &= a + b\tau_1, \\
\tau_3 &= a' + b'\tau_1, \\
\tau_2^2 - \tau_1 \tau_3 &= a'' + b'' \tau_1.
\end{align*}

On the other hand,
\[
 a'' + b'' \tau_1 = \tau_2^2 - \tau_1 \tau_3 = (b^2 - b')\tau_1^2 + (2ab - a')\tau_1 + a^2.
\]
Then either $\tau_1$ solves a quadratic equation (remark that $\tau_1, \tau_3 \in \H$, so in that case $\tau_1$ is quadratic imaginary), and as $\bt \in \M_2(\Q(\tau_1))$  we are finished by Lemma \ref{LemaSimplestPM}, or $b^2 - b' = 0$. The latter is impossible, as $\Im\bt$ is a positive definite symmetric matrix:
\[
\Im \bt = \begin{pmatrix}
\Im(\tau_1) & b\Im(\tau_1) \\ b\Im(\tau_1) & b'\Im(\tau_1)
\end{pmatrix},
\]
\noindent so in particular it has positive determinant, hence $(b' - b^2)\Im (\tau_1) > 0$, so $b^2 \not= b'.$ 
\end{proof}

This result already suggests that $\L_{\bt}$ is a more suitable object for the $\Sp_4(\Z)$-action. More is true:

\begin{lema} If $\bt, \bt' \in \H_2$ parametrize isomorphic ppas, then $(\L_{\bt}, \Delta)$ and $(\L_{\bt'}, \Delta)$ are isomorphic as positive definite lattices. More generally, if they are isogenous, there is $\Z$-linear map $(\L_{\bt}, \Delta) \to (\L_{\bt'}, \Delta')$ such that $\Delta'(\cdot) = \kappa \Delta (\cdot)$ for some constant $\kappa \geq 0$, which extends to a $\Q$-linear isomorphism $\L_{\bt} \otimes \Q \to \L_{\bt'} \otimes \Q$. In particular, $\rank \L_{\bt}$ is invariant under isogenies. \label{LemaRankInvIsogenies} 
\end{lema}

The proof of Lemma \ref{LemaRankInvIsogenies} will be split into two parts. Let us rewrite the statements in terms of actions on $\H_2$. Consider $\bt \in \H_2$ and $(A_{\bt}, H)$ the corresponding principally polarized abelian surface, and assume that we have an isogeny $\phi:B \to A_{\bt}$ with exponent $e(\phi)$, \textit{i.e.} there exists $\psi: A_{\bt} \to B$ with $\psi \phi = e(\phi)_{B}$ and $\phi \psi = e(\phi)_{A}$, by \cite[Proposition 1.2.6]{Birk-LangeCAV}. Hence, $\phi$ induces a polarization of degree $e(\phi)^4$.

The following comes from the proof of \cite[Proposition 8.1.2]{Birk-LangeCAV}. By the theory of elementary divisors, the induced polarization on $B$ has a type $D = \operatorname{diag}(d_1,d_2)$ with $d_1 \lvert d_2$ and $\det(D) = e(\phi)^4$. We can choose a symplectic basis for $A_{\bt}$ and $B$ such that have matrices $R \in \M_4(\Z)$ and $L \in \M_2(\C)$  (the rational and analytic representation of $\phi$, respectively) and period matrices
\[
L(\bt' \: D) = (\bt \: I_2) R.
\]
If we set
\[
{}^{t}\!M = R \begin{pmatrix}
I_2 & 0 \\ 0 & D^{-1}
\end{pmatrix},
\]
then
\[
L(\bt' \: I_2) = (\bt \: I_2) {}^{t}\!M,
\]
and by the same arguments to prove that isomorphisms of principally polarized abelian varieties correspond to matrices in $\Sp_4(\Z)$ acting on $\H_2$, one proves that $M \in \Sp_4(\Q)$ and $\bt' = M\bt$, via a linear fractional transformation.

We claim that $\rank \L_{\bt}$ is invariant by the action of matrices $M \in \Sp_4(\Q)$ such that there exists $D = \operatorname{diag}(d_1, d_2)$, $d_1 \lvert d_2$ and $d_1 d_2 \in \N$ a perfect square with 
\begin{equation}
{}^{t}\!M \begin{pmatrix}
I_2 & 0 \\ 0 & D
\end{pmatrix} \in \M_4(\Z). \label{EqMatrixOfIsogeny}
\end{equation}
Let us first prove the invariance under the $\Sp_4(\Z)$-action separately.

\begin{lema} For $\bt \in \H_2$ and $\bt'=M\bt$ with $M \in \Sp_4(\Z),$ $(\L_{\bt}, \Delta)$ and $(\L_{\bt'}, \Delta)$ are isomorphic as positive definite lattices.
\label{LemmaRankLbt}
\end{lema}

\begin{proof}
By \eqref{EqRationalRepUnderSp4} and \eqref{EqAnalRepUnderSp4}, the action of $M$ in both the rational and analytic representation corresponds to conjugation by suitable matrices. On the other hand, one can check $f \in \End^s(A_{\bt})$ if and only if $f_{M} := \phi^{-1} f \phi \in  \End^s(A_{\bt'})$ for $\phi: A_{\bt'} \to  A_{\bt}$ the isomorphism corresponding to $M$, by the effect on $\phi$ on the principal polarization and on the Rosati involution (equivalently, one can prove it directly for the rational representations $\rho_{r,\bt}(f)$ and $\rho_{r, \bt'}(f)$, for that the corresponding matrix equation \eqref{EqSymRatRep} is invariant under conjugation by matrices in $\Sp_4(\Z)$). 

For $f \in \End^s(A_{\bt})$, by (the proof of) Lemma \ref{LemaHumertRelations}, we read the Humbert singular relation as the non-trivial equation in 
\[
\bt B \bt - C = \bt A - {}^t\!(\bt A),
\]
where $\rho_{r, \bt}(f) = \begin{pmatrix}
A & B \\ C & {}^t\!A
\end{pmatrix}$.

Therefore, as $f_M \in \End^s(A_{\bt'})$ with rational representation given by \eqref{EqRationalRepUnderSp4}, we have a Humbert singular relation for $\bt'$ given by
\[
\bt' B' \bt' - C' = \bt' A' - {}^t\!(\bt' A'),
\]
with coefficients given by ${}^{t}\!M^{-1} \rho_{r,\bt}(f){}^{t}\!M$. Remark that $\rho_{r, \bt}(f) = nI_4$ for some $n \in \Z$ if and only if $\rho_{r, M\bt}(f) = nI_4$, so $\rank \L_{\bt} =0$ if and only if $\rank \L_{\bt'} =0$. Furthermore, $\Delta(f) = \Delta(f_M)$, because $\Delta(f) = \tr(\rho_{a, \bt}(f))^2 - 4\det\rho_{a, \bt}(f)$ is invariant under conjugation. 

Finally, suppose that we have two endomorphisms $f, g \in \End^s(A_{\bt})$, and consider the corresponding Humbert singular relations $l_f, l_g$. They are linearly dependent if and only if there exists $n,m,l \in \Z$ such that
\begin{equation}
n\rho_{r,\bt}(f) = m\rho_{r,\bt}(g) + lI_4. \label{EqLinearDependenceMatrices}
\end{equation}
This follows from the fact that a Humbert singular relation characterizes the matrix of the rational representation up to addition by a multiple of $I_4$, and from $\rho_{r,\bt}(sf) = s\rho_{r,\bt}(f)$ for all $s\in \Z.$ As \eqref{EqLinearDependenceMatrices} is invariant under conjugation, $f, g \in \End^s(A_{\bt})$ produce linearly independent Humbert singular relations if and only if the same is true for $f_M, g_M \in \End^s(A_{\bt'})$. It then follows $\rank \L_{\bt} = \rank  \L_{\bt'}$.
\end{proof}

\begin{proof}[Proof of Lemma \ref{LemaRankInvIsogenies}] Consider $M$ as in \eqref{EqMatrixOfIsogeny} and the set-up described after Lemma \ref{LemaRankInvIsogenies}. For an endomorphism $f \in \End^s(A_{\bt})$ with rational representation $\rho_{r,\bt}(f)$, could consider $\rho_{r,\bt'}(f)$, which gives the rational representation of $f$ in the basis $(\bt' \: I_2)$, \textit{i.e.} $\rho_{r,\bt'}(f) = {}^{t}\!M^{-1} \rho_{r,\bt}(f) M$. As $M \in \Sp_4(\Q)$, it is symmetric with respect to the (extension of) the Rosati involution to $\End^0(A_{\bt'})$, but it is only defined on $\End^0(A_{\bt'})$, as we have formally inverted the isogeny $\phi$. However, we just need to "clear denominators" and consider the correct multiple of $f$ so that we get an honest endomorphism in $\End^s(A_{\bt'})$.

By setting \[
R = {}^{t}\!M \begin{pmatrix}
I_2 & 0 \\ 0 & D
\end{pmatrix} \in \M_4(\Z),
\]
we recover the rational representation of the isogeny (with integer coefficients), hence 
\begin{align*}
\rho_{r,\bt'}(f) = \begin{pmatrix}
I_2 & 0 \\ 0 & D
\end{pmatrix} R^{-1} \rho_{r,\bt} (f) R \begin{pmatrix}
I_2 & 0 \\ 0 & D
\end{pmatrix}^{-1}.
\end{align*}
It is easy to check that although $R^{-1} \not \in \M_4(\Z)$, we have $\begin{pmatrix}
D & 0 \\ 0 & D
\end{pmatrix} R^{-1} \in \M_4(\Z).$

In short, multiplying by a power of $e(\phi)$ will make our matrix integral. From another point of view, in our chosen basis $\rho_{r,\bt'}(f)$ is the rational representation of $\psi^{-1} f \psi$, formally inverting $\phi$, so we need to instead consider $\phi f \psi$ (from this description we can also see that $\phi f \psi$ is symmetric). Because $\psi^{-1} = \frac{1}{e(\phi)}\phi$, this means that 
\[
e(\phi) \rho_{r,\bt'}(f) \in \M_4(\Z).
\]

Now that we have defined a map $\L_{\bt} \to \L_{\bt'}$, we can argue, as in the proof of Lemma \ref{LemmaRankLbt}, that $\rank \L_{\bt'} = \rank \L_{\bt}$, and $\Delta' = e(\phi)^2 \Delta$.

As an alternative proof of the invariance of the rank under isogenies, note that by Proposition \ref{PropIsomNSEnds}, there exists an isomorphism of abelian group between $\NS(A_{\bt})$ and $\End^s(A_{\bt})$. In addition, by \eqref{EqSESNeronSeveri}, $\rank \L_{\bt} = \rank \NS(A_{\bt}) - 1$. Finally, it is known that the rank of the N\'eron-Severi group of an abelian variety is invariant under isogenies, see \cite[Chapter 1, Prop 3.2]{BirkLangeCT}.
\end{proof}

We devote the rest of this section to the computation of $\rank(\L_{\bt})$. 

\begin{prop} If $\bt \in \H_2$ corresponds to a ppas $A_{\bt}$ such that $\End_0(A_{\bt})$ is commutative but distinct from $\Q$, which means that it is either
\begin{enumerate}
\item a real quadratic field,
\item a CM quartic field,
\item $\Q \times \Q$,
\item $\Q \times K$ for $K$ an imaginary quadratic,
\item or $K_1 \times K_2$ for $K_1 \not= K_2$ imaginary quadratic fields,
\end{enumerate} 
then $\rank \L_{\bt} = 1.$\label{PropRank1Case}
\end{prop}

\begin{proof}
We remark that $\rank \L_{\bt} = 1$ is equivalent to $A_{\bt}$ solving a unique (up to sign) primitive HSR, hence $A_{\bt}$ belonging to a unique Humbert surface of minimal discriminant.

The cases (1) and (2) correspond to $\End(A)$ without zero divisors (equivalently $A_{\bt}$ simple), then there cannot exist $f \in \End_0(A_{\bt})$ such that $\Delta(f)$ is a perfect square, by Proposition \ref{PropDeltasAndQuadFieldIsogenies}. In addition, there is a \emph{unique} quadratic real field embedded in $\End_0(A_{\bt})$: in the case (1) $\End_0(A_{\bt})$ is already a quadratic field, and in (2) it corresponds to the (unique) maximal real subfield of the CM field. By Albert's classification of endomorphisms algebras of simple abelian varieties \cite[Proposition 5.5.7]{Birk-LangeCAV}, in (1) the Rosati involution is trivial, and in (2) the Rosati involution corresponds to complex conjugation, so in both cases the embeddings restrict to $\End_0^s(A_{\bt})$, and are surjective. Hence $\End_0^s(A_{\bt})$ is a real quadratic field, and $\End^s(A_{\bt})$ is an order in said field. Taking $f \in \End^s(A_{\bt})$ with minimal non-zero $\Delta(f)$ (the discriminant of the order by which $A_{\bt}$ has real multiplication), we necessarily have $\{n + mf\} = \End^s(A)$. Therefore, we have one (up to sign) primitive HSR, and any other is necessarily a multiple of it.

The other cases correspond to $A_{\bt} \simeq E \times E'$ with $E$ and $E'$ nonisogenous elliptic curves. First, $\End^s_0(A_{\bt})$ has (non trivial) symmetric idempotents by \cite[Theorem 5.3.2]{Birk-LangeCAV}, in particular there exists $f \in \End^s(A_{\bt})$ with $f \not\in \Z$, so that $\rank\L_{\bt} \geq 1$.  We note that $A_{\bt}$ cannot have real multiplication by any real quadratic field, as by \cite[Corollary 2.7]{Goren}, the abelian surfaces with real multiplication by a quadratic field are either simple or isotypic (isogenous to the square of an elliptic curve).

Therefore, $\Delta$ induces a positive definite quadratic form on $\L_{\bt}$ that only represents squares. By contradiction, assume that there is a sublattice of $\L_{\bt}$ of rank two, then $\Delta$ will induce a positive definite \emph{binary} quadratic form $rf(x,y)$ with $f = ax^2 + bxy + cy^2$ primitive (\textit{i.e.} $\gcd(a,b,c) = 1$) and $r \in \N$. It is known that primitive positive definite quadratic forms represent infinitely many prime numbers (see \cite[Theorem 9.12]{Cox}), which is a contradiction to $rf$ only representing squares. Therefore, $\rank \L_{\bt} = 1.$

\end{proof}

The lattice of singular relations is a more interesting object in the quaternionic multiplication case. Note that if we assume $A$ ppas admitting an embedding $B \hookrightarrow \End_0(A_{\bt})$, then either $B \cong \End_0(A)$ or $\End_0(A)$ is strictly larger (and in particular, $\rank_\Q (\End_0(A)) \geq 4$). In the latter case, one can check by inspection with all the other possibilities for $\End_0(A)$ with $\rank_\Q (\End_0(A)) \geq 4$, that necessarily $\End_0(A) \cong M_2(K)$ for $K$ a imaginary quadratic field. Therefore $A_{\bt} \simeq E^2$ for $E$ an elliptic curve with CM by $K$, and by Lemma \ref{LemaRankLtExtremeCases}, it is completely characterized by $\rank \L_{\bt} = 3$.

\begin{prop} Let $B$ be an indefinite quaternion \emph{division} algebra over $\Q$. Suppose that $A_{\bt}$ admits quaternionic multiplication by an order in $B$.
\begin{itemize}
\item If $A_{\bt}$ is simple (equivalently $\End_0(A) \cong B$) then $\rank \L_{\bt} = 2$.
\item Otherwise, $A_{\bt}$ is CM and $\rank \L_{\bt} = 3$.
\end{itemize} \label{PropQMcase}
\end{prop}
\begin{proof}
If $A_{\bt}$ is simple, by \cite[Lemma 7.9]{RemVarAbelOrdMax}, then we can fix some Eichler order $\O$ (intersection of maximal orders) of square-free level in $B$, and there exists an isogeny from $A_{\bt}$ to an abelian surface $A'$ with $\End(B) = \O$. This abelian variety $A'$ further admits a principal polarization, as we will see in Section \ref{sec:SecShimurCurves}. By Lemma \ref{LemmaRankLbt}, it is enough then to compute $\rank \L_{\tilde{\bt}}$, with $\tilde{\bt}$ a period matrix for $A'$. In Section \ref{sec:SecShimurCurves}, we will state more information about ppas with quaternionic multiplication by Eichler orders; in particular, $\tilde{\bt}$ is in the $\Sp_4(\Z)$-orbit of a matrix in the image of the morphism in Proposition \ref{PropHas95}, from which it follows that $\rank \L_{\tilde{\bt}} = 2$. By Lemma \ref{LemaRankInvIsogenies}, $\rank \L_{\bt} = 2.$

Alternatively, for a more intrinsic argument, we could have also fixed $\O$ as a maximal order in $B$. In the case of maximal orders, by \cite[Lemma 16 and Lemma 17]{Lin-Yang}, the lattice of singular relations is isomorphic to a lattice of rank two contained in $B$.

The second statement was proven in Lemma \ref{LemaRankLtExtremeCases}.
\end{proof}

To finish the study of $\mathcal{L}_{\bt}$, we are only missing the isotypic case $A_{\bt} \simeq E^2$, where $E$ is \emph{without} complex multiplication. By Lemma \ref{LemaRankInvIsogenies}, assume that $A_{\bt}$ is \emph{isomorphic} to the product of elliptic curves with product polarization.

\begin{lema} Consider an elliptic curve $E$ and the abelian surface $E^2$ (with the canonical principal polarization). There exists an embedding
\[
\mathcal{O}_{K} \hookrightarrow \End(E^2),
\]
for all real quadratic fields $K$. \label{LemaE2RMByAllQuadFields}

With respect to the standard principal polarization, there exist infinitely many quadratic fields such that 
\[
\mathcal{O}_{K} \hookrightarrow \End^s(E^2).
\]
\end{lema}

\begin{rem} The endomorphism algebra of an abelian variety is invariant under isogenies, but that is not true for the endomorphism ring. What we can say is that if $A_{\bt} \simeq E^2$, then $A_{\bt}$ admits real multiplication by an \emph{order} in $K$ for any real quadratic field $K$, but not necessarily by the ring of integers $\O_K$.
\end{rem}

\begin{proof}
This is a classical observation from the theory of real multiplication, for example in \cite[Example 2.2 (2)]{Goren}, if $d$ is square free, we can choose integers $a,b,c$ with $d = a^2 + bc$. 

If $d \equiv 2,3 (4)$ then the embedding from $\O_{\Q(\sqrt{D})} \to \M_2(\Z)$ is given by 
\[
\sqrt{d} \mapsto \begin{pmatrix}
a & b \\ c & -a
\end{pmatrix}.
\]

If $d \equiv 1 (4)$ we can choose $a$ odd and $b,c$ even (take $a=1$, $b=2$, $c = (D-1)/2$), and the embedding corresponds to 
\[
\frac{1 + \sqrt{d}}{2} \mapsto \begin{pmatrix}
\frac{1 + a}{2} & \frac{b}{2} \\ \frac{c}{2} & \frac{1-a}{2}
\end{pmatrix}.
\]

The canonical polarization induces the Rosati involution $\M_2(\Q) \to \M_2(\Q)$ given by $Z \mapsto {}^t\!Z$. Hence, the previous embeddings map to $\End^s(A_{\bt})$ when can take $b=c$ in $d = a^2 + bc$. By Jacobi's two square theorem, it follows that $\O_K \hookrightarrow \End^s(A_{\bt})$ for square free integers $d = 2^jp_i \cdots p_r$ with $j=0,1$ and all $p_i \equiv 1$ (mod $4$).

\end{proof}
\begin{prop} If $A_{\bt} \simeq E^2$, where $E$ does \emph{not} have CM, then $\rank \mathcal{L}_{\bt} = 2.$ \label{PropIsogCaseRank2}
\end{prop}

\begin{proof}
As we have said previously, it is enough to show it for $E^2$ with the canonical product polarization. By Lemma \ref{LemaRankLtExtremeCases}, it cannot be $0$ or $3$. Assume that $\rank \mathcal{L_{\bt}} = 1$, hence there exists a primitive singular relation $l$ such that $\L_{\bt} = ml$ for $m \in \Z$, and $\Delta$ restricts to the quadratic form $f(m) = m^2\Delta(l)$, where $\Delta(l)$ is a fixed constant. On the other hand, by Lemma \ref{LemaE2RMByAllQuadFields}, $m^2\Delta(l)$ must represent infinitely many fundamental discriminants, which is a contradiction. Therefore, $\rank \mathcal{L}_{\bt} = 2.$
\end{proof}

We collect all the past results (Lemma \ref{LemaRankLtExtremeCases}, Proposition \ref{PropRank1Case}, Proposition \ref{PropQMcase} and Proposition \ref{PropIsogCaseRank2}) of this section in the following theorem. Remark that the proofs of 2) and 3) are for one implication only, but as we have exhausted all the possible cases for $\End_0(A_{\bt})$, they are also equivalences.
\begin{teo} Let $\bt \in \H_2$ and consider the corresponding ppas $A_{\bt} \in \A_2$.
\begin{itemize}
\item It holds $\End_0(A_{\bt}) = \Q$ if and only if $\rank \mathcal{L}_{\bt} = 0$.
\item It holds $\End_0(A_{\bt})$ is commutative but not $\Q$ if and only if $\rank \mathcal{L}_{\bt} = 1.$
\item It holds $\End_0(A_{\bt})$ is an indefinite quaternion algebra over $\Q$ if and only if $\rank \mathcal{L}_{\bt} = 2.$
\item It holds $\End_0(A_{\bt}) = \M_2(K)$, for $K$ a quadratic imaginary field, if and only if $\rank \mathcal{L}_{\bt} = 3.$
\end{itemize} \label{TeoBigThingHSR}
\end{teo}

Alternatively, this is also proven in \cite[Chapter 2, Proposition 7.1]{BirkLangeCT}, by computing $\rank \NS(A)$. We underline the case of $\End_0(A_{\bt})$ with zero divisors, which follows from a computation of $\rank \NS(E\times E')$ with $E,E'$ elliptic curves. By \cite[Proposition 23]{KaniThemodulispace} there is a group isomorphism:
\[
\NS(E \times E') \cong \Z \oplus \Z \oplus \Hom(E,E').
\]
Therefore, for $A_{\bt} \simeq E \times E'$ (with $E$ and $E'$ not necessarily isogenous) depending on \\$\rank \Hom(E,E') = 0,1,2$, we have
\[
\rank \L_{\bt} = \begin{cases} 1, & \text{ if } E \not\simeq E', \\
2, & \text{ if } E \simeq E', \text{ without CM}, \\
3, & \text{ if } E \simeq E', \text{ with CM}.
\end{cases}
\]

\section{From Humbert singular relations to linear relations}
We now have a good understanding of $\L_{\bt}$ from Theorem \ref{TeoBigThingHSR}. However, we were initially interested in $\Llin_{\bt}$. This section is devoted to infer the relevant information about it to prove Theorem \ref{TeoBigTeo}.

Let us first notice that $\rank \Llin_{\bt} \leq \min(\rank \L_{\bt}, 2)$, and if $\rank \L_{\bt} = 0$, then $\rank \Llin_{\bt} = 0$ too, and since $\rank \L_{\bt}$ is invariant under $\Sp_4(\Z)$-orbit, then $\rank \Llin_{\bt} = 0$ holds for the whole orbit of $\bt$.

\subsection{The case \texorpdfstring{$\rank \L_{\bt} = 1$}{rank L = 1}, and proof of Theorem \ref{TeoBigTeo} 1), 2)}

If $\bt \in \H_2$ is a CM point and $\rank \L_{\bt} = 1$ then there exists (up to sign) only one primitive Humbert singular relation $l$ satisfied by $\bt$. It can involve $\t_2^2 - \t_1\t_3$ or not. In the latter case, $\dim \spbt = 3$ and by Corollary \ref{CorConjeForCM}, $\trqjq = 2$.

If $l$ involves the determinant, then by Humbert's lemma Proposition \ref{PropHumbertsLemma}, there exists $\bt'$ in the $\Sp_4(\Z)$-orbit with $\dim \spn_\Q(1, \t'_1, \t'_2, \t'_3) = 3$, and can apply the previous argument now to $\bt'$. It then follows that $\min_{\Sp_4(\Z)} \trqjq = 2.$

Set $\Delta = \Delta(l)$ and consider the Humbert surface $\Hum_{\Delta}$. Consider now other $\tilde{\bt} \in \Hum_{\Delta}$ (non-CM) such that $j_1(\tilde{\bt}), j_2(\tilde{\bt}) , j_3(\tilde{\bt})  \in \Qbar$. Then $\trqjqtil = \trdeg \Q(\tilde{q}_1, \tilde{q}_2, \tilde{q}_3)$. On the other hand, they also satisfy the Humbert singular relation $l$, so we can argue as above and have $\min_{\Sp_4(\Z)} \trqjqtil \leq 2$.

We have therefore proven the first two statements in Theorem \ref{TeoBigTeo}. In the second statement, note that if $A_{\bt}$ belongs to a special curve of type $\Q \times CM$, this special curve can also be taken as the special subvariety, for that $\rank \L_{\bt} = 1$ generically in that curve by Proposition \ref{PropRank1Case} (4). 

\subsection{The case \texorpdfstring{$\rank \L_{\bt} = 2,3$}{rank L = 2,3}} \label{sec:SecRankHigh}

If $\rank \L_{\bt} = 3$, then necessarily $\rank \Llin_{\bt} = 2.$ Also, by Theorem \ref{TeoBigThingHSR}, this case can only happen for $A_{\bt}$ isotypic CM.

If $\rank \L_{\bt} = 2$ then $\rank \Llin_{\bt} = 2$ or $1$, and we will show that both cases occur. By Proposition \ref{PropRank1Case}, they have to be either Shimura curves or modular curves. For the proof of Theorem \ref{TeoBigTeo}, we want a special curve such that $\rank \Llin_{\bt}=2$, so our discussion in Section \ref{sec:SecShimurCurves} is not relevant for the proof of Theorem \ref{TeoBigTeo}, though it is still insightful to complete the study for $\Llin_{\bt}$.

\subsubsection{Shimura curves by an Eichler order} \label{sec:SecShimurCurves}

For the rest of the section, we need to introduce the quaternionic modular embeddings for an Eichler order in an indefinite quaternion algebra over $\Q$. This construction is standard in the literature for Shimura curves, see, for example, \cite[Section 2.1]{Guo-Yang} or \cite[Section 3]{Hashimoto}, and comes from a classical paper by Shimura \cite{ShimuraAnalyticFamilies}.

Let $B$ an indefinite quaternion algebra over $\Q$ of discriminant $D$ (with $D > 1$, equivalently, $D$ is a division algebra). Let $\mathcal{O}$ be an Eichler order of level $N$ squarefree (by definition, an intersection of two maximal orders). 

Set for $\tau \in \H$ the vector $v_{\tau} = (\tau \; 1)^t \in \C^2$, and fix an embedding $\phi: B \to \M_2(\R)$. We remark that such an embedding exists as $B$ is indefinite, \textit{i.e.} $B \otimes_{\Q} \R \cong \M_2(\R)$, and any two embeddings differ by conjugation by an invertible matrix in $\M_2(\R)$, by the Skolem-Noether theorem.

We want to assign to every $\tau \in \H$ a complex torus with a principal polarization. We set:
\begin{itemize}
\item the lattice $\Lambda_{\tau} := \phi(\mathcal{O}) v_{\tau}$, with $\phi(\mathcal{O})\subset \M_2(\R)$ acting as $2 \times 2$ matrices, and the torus $\faktor{\C^2}{\Lambda_{\tau}}$;
\item for the principal polarization, let $\mu \in \mathcal{O}$ satisfying the following conditions
\begin{itemize}
\item it holds $\tr(\mu) = 0$, $\nrd(\mu) > 0$, and $\mu^2 + DN = 0$, and
\item if we set $\phi(\mu) = \begin{pmatrix}
a & b \\ c & -a
\end{pmatrix}$ then $c > 0$,
\end{itemize} 
then there is a symplectic form $E_{\mu}$ on $\Lambda_{\tau} \times \Lambda_{\tau} \to \Z$ (and extended to $\C^2$ by $\R$-linearity) defined by
\[
E_{\mu}(\phi(\alpha) v_{\tau}, \phi(\beta) v_{\tau}) := \tr(\mu^{-1}\alpha \bar{\beta}).
\]
The conditions imposed to $\mu$ make $E_{\mu}$ a Riemann form, and $A_{\tau} = (\C^2/ \Lambda_{\tau}, E_{\mu})$ a principally polarized abelian surface. Likewise, we have an embedding $\O \hookrightarrow \End(A_{\tau})$ induced by $\phi$, compatible with the Rosati involution induced by $E_{\mu}$. The restriction of the Rosati involution to $\O$ is given by $\al \mapsto \mu^{-1} \bar{\al} \mu$. 
\end{itemize}

This assignment induces a well-defined map $X_0^{D}(N) \to \A_2$, for $X_0^{D}(N)$ the Shimura curve associated to $\O$. Set $\mathfrak{X}_{\mu}$ for its image in $\A_2$.

We can further define a map $\H \to \H_2$ by choosing a symplectic basis $\al_1, \ldots , \al_4$ in $\mathcal{O}$ with respect to $E_{\mu}$. Then the big period matrix is given by $\Pi(\tau) = (\al_1 v_{\tau}, \al_2 v_{\tau}, \al_3 v_{\tau}, \al_4 v_{\tau}) = \left( \Omega_1 (\tau), \Omega_2(\tau)\right)$, and the period matrix by $ \Omega_2(\tau)^{-1} \Omega_1(\tau).$ This basis also induces an embedding via $\phi$ of $\O^1$ into $\operatorname{Aut}_{\Z}(\O, E_\mu) \cong \Sp_4(\Z)$, see \cite[Proposition 2.5]{Hashimoto}.

Those maps can be made explicit, see \cite[Theorem 3.5 and Theorem 5.1]{Hashimoto}. It then follows that $\rank \L_{\bt} = 2$ and $\rank \Llin_{\bt} = 1$ for $\bt$ in the image of these embeddings. There is a choice of an odd prime $p$ and an integer $a$, see \cite[page 535 above Equation (1)]{Hashimoto} for more details. Likewise, see \cite[Lemma 38]{Guo-Yang} and \cite[Lemma 5]{Lin-Yang}, where these choices are reinterpreted in terms of $\GL_2(\Z)$-equivalence classes of positive definite binary quadratic forms of discriminant $-16DN$.

\begin{prop}[\cite{Hashimoto}] Let $\mathcal{O}$ and $B$ be as above, and set 
$\eps = \frac{1 + \sqrt{p}}{2}$ and $\bar{\eps} = \frac{1 - \sqrt{p}}{2}$. Then the following map $\Omega$ gives a modular embedding of $\H$ into $\H_2$ with respect to $\phi(\O^1)$ and $\Sp_4(\Z):$ 
$$
\Omega(z) = \frac{1}{p z}\left(\begin{array}{cc}
-\bar{\varepsilon}^2+\frac{(p-1) a DN}{2} z+DN \varepsilon^2 z^2, & \bar{\varepsilon}-(p-1) a DN z-DN \varepsilon z^2 \\
\bar{\varepsilon}-(p-1) a DN z-DN \varepsilon z^2, & -1-2 a DN z+DN z^2
\end{array}\right).
$$  
Furthermore, $\bt=\Omega(z)$ satisfy simultaneously the following singular relations parametrized by two independent integers $x, y \in \Z:$
$$
x \tau_1+(x+2 a DN y) \tau_{2}-\frac{p-1}{4} x \tau_3+y\left(\tau_{2}^2-\tau_1 \tau_3\right)+\left(a^2 DN-b\right) DN y=0,
$$
where we put $a^2 DN+1=p b$. Moreover, if $z \in \H$ is not a CM point, then it does not have another singular relation. \label{PropHas95}
\end{prop}

One can then see that $\rank {\L^{lin}_{\bt}} = 1$ for $\bt$ in $\mathfrak{X}_{\mu}$.

Let us understand the subset covered by $\cup_{\mu}\mathfrak{X}_{\mu}.$ The following is a consequence of \cite[Section 2, Theorem 1 and Theorem 2]{ShimuraAnalyticFamilies}, that we state as in \cite[page 6, paragraph before section 2.2]{Guo-Yang}. Set $\mathcal{Q}_{D,N} \subset \A_2$ the set of ppas with QM by an Eichler order $\O$ of level $N$ in a quaternion algebra of discriminant $D$, such that the Rosati involution restricted to $\O$ coincides with $\al \mapsto \mu^{-1}\bar{\al}\mu$, for $\mu \in \O$ allowed as above \footnote{As by \cite[\S 43.6.2]{VoightQuatAlg}, when $\End(A) = \O$, or equivalently for the simple abelian surfaces in $\mathcal{Q}_{D,N}$, the last condition is always satisfied.}. Then
\begin{equation}
\mathcal{Q}_{D,N} = \bigcup_{\mu \in \O}  \mathfrak{X}_{\mu}, \label{EqQuaternLocus}
\end{equation}
where $\mu$ runs through the allowed elements in $\O.$ 

We can have two distinct allowed elements $\mu_1 \not= \mu_2 \in \O$ such that $\mathfrak{X}_{\mu_1} = \mathfrak{X}_{\mu_2}.$ Actually, there are only \emph{finitely many} connected components in \eqref{EqQuaternLocus}, by \cite[Proposition 4.3]{RotgerShimCurIgusa}.

\begin{rem} As $\rank \Llin_{\bt}$ is not invariant under the $\Sp_4(\Z)$-orbit, we can only say that for this explicit period matrix, $\rank \Llin_{\bt} = 1$, but it could be possible that for another element of the orbit $\rank \Llin_{\bt} = 2$. 

We conjecture that this is not the case, and that $\rank \Llin \not= \rank \L$ throughout the $\Sp_4(\Z)$-orbit. We partially solve this for some Shimura curves in the Appendix. \label{RemConjNonSimHumbForShimCurves}
\end{rem} 

\subsubsection{A collection of modular curves}  
\label{sec:SecFamily}

The modular embedding construction for Shimura curves $X_0^D(N)$ can be applied to $\M_2(\Q)$ and $\mathcal{O}_{N} := \begin{pmatrix}
\Z & \Z \\ N\Z & \Z
\end{pmatrix}\subset \M_2(\Q)$, $N$ squarefree. We follow here \cite[Section 6]{Lin-Yang} for the case $N$ squarefree, and \cite{KaniThemodulispace} for general $N$.

We take $\mu \in \mathcal{O}_{N} := \begin{pmatrix}
\Z & \Z \\ N\Z & \Z
\end{pmatrix}\subset \M_2(\Q)$ with trace $0$, determinant $N$ and positive $(2,1)$ entry, \textit{i.e.}
\[
\mu = \begin{pmatrix}
a & b \\ cN & -a
\end{pmatrix},
\]
with $c > 0$ and $-a^2 - bcN = N$. This implies $b<0$ and, as $N$ is squarefree, we have $N \mid a$. Let us rewrite $\mu$ as 
\[
\mu = \begin{pmatrix}
aN & -b \\ cN & -aN
\end{pmatrix},
\]
with $b,c >0$ and $a\in \Z$ such that $bc - Na^2=1$.
\begin{lema} For $\mu \in \O_N \subset \M_2(\Q)$ as above, we have a map $Y_0(N) \to \A_2$ induced from
\begin{align*}
\H &\to \H_2 \\
\tau &\mapsto \begin{pmatrix}
b\tau & aN\tau \\ aN\tau & cN\tau
\end{pmatrix}.
\end{align*}
For $\bt$ of this shape, we have $\rank \L_{\bt} = \rank \Llin_{\bt} = 2.$ \label{LemaFamilyModCurves}
\end{lema}

\begin{rem} Taking $\mu = \begin{pmatrix}
0 & -1 \\ N & 0
\end{pmatrix}$ we recover the diagonal embedding
\[
\tau \mapsto \begin{pmatrix}
\tau & 0 \\ 0 & N\tau 
\end{pmatrix}
\]
of $Y_0(N) \to \A_1 \times \A_1$, with the product of elliptic curves and product polarization. However, different choices of $\mu$ allow us to consider modular curves contained in the indecomposable locus, so we can evaluate the Igusa invariants on them.

If $N$ is not square-free, Lemma \ref{LemaFamilyModCurves} still constructs some explicit maps from $Y_0(N)$ to $\A_2$, although not necessarily all of them. 

Finally, notice that the explicit map in Lemma \ref{LemaFamilyModCurves} gives an alternative proof of Corollary \ref{CorCompEisotypic}, for the images of quadratic imaginary $\tau \in \H$. \label{RemDegenerateModularLocus}
\end{rem} 
\begin{proof}
Set the embedding $\phi: \M_2(\Q) \hookrightarrow \M_2(\R)$ as simply the natural inclusion. We remark that we are considering complex tori $\C^2/ \Lambda{\tau}$ where $\Lambda{\tau} = \mathcal{O}_{N}v_{\tau}$, and the Riemann form $E_{\mu}(\alpha v_{\tau}, \beta v_{\tau}) = \tr(\mu^{-1}\alpha \bar{\beta})$. Note that the quaternionic conjugation in $\M_2(\Q)$ is adjugation of matrices.

One can check by computation that the following basis for $\mathcal{O}_{N}$
\[
\al_1 = \begin{pmatrix}
b & 0 \\ aN & 0
\end{pmatrix}, \:  \al_2 = \begin{pmatrix}
aN & 0 \\ cN & 0
\end{pmatrix}, \:  \al_3 = \begin{pmatrix}
0 & 1 \\ 0 & 0
\end{pmatrix}, \:  \al_4 = \begin{pmatrix}
0 & 0 \\ 0 & 1
\end{pmatrix},
\]
is a symplectic basis with respect to $E_{\mu}.$ It is a basis because from $bc - Na^2=1$ it follows that $\gcd(b,aN) = \gcd(a,c) = 1$, and we have verified that this basis is symplectic with Sage. Here is the code
\begin{verbatim}
var('x,y,a,b,c,N')
mu = matrix([[a*N,-b], [c*N,-a*N]])
al1 = matrix([[b, 0],[a*N,0]])
al2 = matrix([[a*N,0],[c*N,0]])
al3 = matrix([[0,1],[0,0]])
al4 = matrix([[0,0],[0,1]])
al = [al1, al2, al3, al4]
def E(x,y):
    return ((mu.inverse())*x*(y.adjugate())).trace()
List = [] #computes the values of E(x,y) to check they are either 0 or 1
List2 = [] #returns the pairs of indices with E(x,y) = 1
for i in range(4):
    List += [E(al[i], al[i]).full_simplify()]
    for j in range(i,4):
        List += [E(al[i], al[j]).full_simplify()]
        if E(al[i], al[j]) != 0:
            List2 += [(i,j)]
List
\end{verbatim}

The big period matrix is then $\Pi(\tau) = (\al_1 v_{\tau}, \al_2 v_{\tau}, \al_3 v_{\tau}, \al_4 v_{\tau})$ and the embedding into $\H_2$ is given by the statement. This choice of symplectic basis also induces a compatible embedding from $\O_N^1 = \Gamma_0(N)$ to $\Sp_4(\Z)$, so it induces a map $Y_0(N) \to \A_2$. It follows that $\det(\bt) = N\tau^2$, and there cannot be a generic linear dependence relation between $\det(\bt)$ and $\t_1, \t_2, \t_3$. 
\end{proof}

We can also explicitly determine the quadratic form associated to $\L_{\bt}$. Solving for the singular relations in $(p,q,r,s,t)$
\[
\left( bp + aNq + cNr \right)\tau + sN\tau^2 + t = 0
\]
forces $s=t=0$ and 
\[
bp = -N\left( aq + cr \right).
\]
As $bc - Na^2 = 1$, in particular $(b,N) =1$, therefore, necessarily $N \lvert p$. One can check that
\[
(0, -c, a, 0 ,0)
\]
and 
\[
(-N, -Nab , b^2,0,0)
\]
are generators of the lattice (every possible value of $p$ is attained, and for a fixed value of $p$, $(q,r)$ are the solutions to a B\'ezout type equation and one can check that all such solutions are generated). Therefore,
\[
\L_{\bt}= \{ (-Ny, -cx -Naby, ax + b^2 y, 0, 0),\; x,y \in \Z \},
\]
and 
\begin{equation}
\Delta(x,y) = c^2 x^2 + 2aN(bc+2)xy + b^2N(bc + 3)y^2. \label{EqStandardQuadraticForm}
\end{equation}

We have recovered the same positive definite quadratic form as in \cite[Equation (5)]{KaniThemodulispace}, with parameters $(c,b,a) \in P(N)$ in their notation. By \cite[Theorem 13]{KaniThemodulispace}, this quadratic form is characterized by the following properties ("quadratic forms of type $N$"):
\begin{itemize}
\item its discriminant as a binary quadratic form is $-16N$,
\item for any $x,y \in \Z$, $\Delta(x,y) \equiv 0,1$ (mod $4$),
\item $\Delta$ primitively represents a square prime to $N$ (note that $\Delta (1,0) = c^2$, $\Delta(b^2, -a) = b^2$).
\end{itemize}

In addition, every such quadratic form of type $N$ is $\SL_2(\Z)$-equivalent to one as in \eqref{EqStandardQuadraticForm} for some parameters $(c,b,a)$ such that $bc - Na^2 = 1.$ We call the quadratic forms in \eqref{EqStandardQuadraticForm} \emph{standard}.

Consider now a general $N$. As we mentioned in Remark \ref{RemDegenerateModularLocus}, not all admissible $\mu$ are covered by the explicit formula in Lemma \ref{LemaFamilyModCurves}. However, it follows from the treatment in \cite{KaniThemodulispace} (where $N$ is not required to be square-free), that it is enough to construct the standard ones, as they serve as a full set of representatives of the finitely many distinct images of these maps, in bijection with the $\GL_2(\Z)$-equivalence class of the quadratic forms of type $N$. 

Analogously to the Shimura curves, let us study the loci $\cup_{\mu} \mathfrak{X}_{\mu}.$ In \cite{KaniThemodulispace}, these loci are considered in $\M_2 \cong \A_2^{ind}$, which is sufficient to us, as we eventually want to evaluate the Igusa invariants. The problem studied is the following: it is classical (\cite[Satz 2]{Weil}) that for $C \in \M_2$, the ppas $(\Jac(C), \theta)$ belongs to the indecomposable locus of $\A_2$, in other words, for every $E,E'$ pair of elliptic curves:
\[
(\Jac(C), \theta) \not\cong (E \times E', L_E \oplus L_{E'}), 
\]
as ppas, where $L_E \oplus L_E'$ is the canonical induced principal polarization on $E \times E'$. However, one could ask if there exist another principal polarization $\tilde{L}$ on $E \times E'$ such that
\begin{equation}
(\Jac(C), \theta) \cong (E \times E', \tilde{L}),\label{EqTorelliProperty}
\end{equation}
or equivalently, if $\Jac(C)$ is isomorphic to $E \times E'$ as \emph{unpolarized} abelian varieties. 

It follows that if \eqref{EqTorelliProperty} holds for some $C \in \M_2$, then necessarily $\Hom(E,E') \not= 0$ by \cite[Proposition 26]{KaniThemodulispace}, or also \cite[Lemma 2.2]{LangePrinPolProdEC} or \cite[Lemma 2.1]{Daw-OrrUI}. Furthermore, if $\Jac(C)$ is not CM, then the elliptic curves are not CM and there exists a \emph{unique} cyclic isogeny of minimal degree $N$ (the degree of an isogeny generating $\Hom(E,E') \cong \Z)$, see \cite[Corollary 27]{KaniThemodulispace}.

\begin{defi}[\cite{KaniThemodulispace}] Let $N \geq 1$ an integer, a curve $C \in \M_2$ \emph{has type} $N$ if there exists $E,E'$ elliptic curves, a cyclic isogeny $\al: E \to E'$ with $N = \deg(\al)$ and an isomorphism of abelian varieties $\Jac(C) \cong E \times E'$. Equivalently, there exists a principal polarization $\tilde{L} = \tilde{L}(\al)$ on $E \times E'$ such that \eqref{EqTorelliProperty} holds for $(\Jac(C), \theta)$. The polarization $\tilde{L}$ depends on $\al$ in an explicit way.
\end{defi}

Consider the locus of $\M_2$, for $N \in \N$  defined by
\[
K(N) = \{ C \in \M_2| \; C \text{ has type } N \}.
\]
By \cite[Theorem 12]{KaniThemodulispace}, these loci $K(N)$ are individually covered by modular curves associated to $\GL_2(\Z)$-equivalence classes of quadratic forms of type $N$, and any such quadratic form is $\SL_2(\Z)$-equivalent to one of the form of \eqref{EqStandardQuadraticForm}. 

One could additionally ask if \emph{all} modular curves are of this shape. We do not need to ponder about that, because it is true that every possible CM point belongs to one of $K(N)$, as a consequence of the following result.

\begin{teo}[Shioda and Mitani]If an abelian surface is isogenous to $E^2$ with $E$ elliptic curve with complex multiplication, then it is \emph{isomorphic} (as unpolarized abelian surfaces) to a product of elliptic curves. \label{TeoShiodaMitani}
\end{teo}
\begin{proof}
See \cite[Corollary 10.6.3]{Birk-LangeCAV}.
\end{proof}

We finish this section with a slightly more general result via an alternative approach. Consider $\bt \in \H_2$ and the lattice $\Lambda = (I_2 \; \bt)\Z^4 \in \C^2$. Remark that $\Lambda \otimes_{\Z} \Q$ admits a natural left $\End_{0}(A_{\bt})$-structure. If $A_{\bt}$ is isotypic, then $\M_2(\Q) \subset \End_{0}(A_{\bt})$, hence $\Lambda \otimes_{\Z} \Q$ admits a $\M_2(\Q)$-structure.

In the previous examples, this $\M_2(\Q)$-structure was given simply by matrix multiplication, as the embedding $\M_2(\Q) \hookrightarrow \M_2(\R)$ was the natural inclusion.

\begin{lema} Let $\bt \in \H_2$ with $A_{\bt} \simeq E^2$. Suppose that the $\M_2(\Q)$-module structure of $\Lambda_{\bt} \otimes \Q$ is given by matrix multiplication. Then $\dim \spn_{\Q}(1, \tau_1, \tau_2, \tau_3) = 2$. \label{LemaGeneralFamilyRank2}
\end{lema}

\begin{proof}
By hypothesis, for any $M \in \M_2(\Q)$, we have $M v \in (I_2 \; \bt)\Q^4$ for $v$ any of the basis vectors. In particular, taking $M = E_{ij}$ elementary matrices, we derive the following relations for $\bt$: there exists $a,b,c,d \in \Q$ such that:
\[
\begin{pmatrix}
\tau_1 \\ 0 
\end{pmatrix} = 
\begin{pmatrix}
1 & 0 \\ 0 & 0
\end{pmatrix} \begin{pmatrix}
\tau_1 \\ \tau_2
\end{pmatrix} = a \begin{pmatrix}
1 \\0 
\end{pmatrix} + b \begin{pmatrix}
0 \\ 1
\end{pmatrix} + c \begin{pmatrix}
\tau_1 \\ \tau_2
\end{pmatrix} + d \begin{pmatrix}
\tau_2 \\ \tau_3
\end{pmatrix}
\] 
In other words,
\begin{align*}
\tau_1 &= a + c\tau_1 + d\tau_2,\\
0 &= b + c\tau_2 + d\tau_3.
\end{align*}
From the second equation, we have $c \not= 0$, because $\tau_3 \not\in \Q$ (and $b=d=c=0$ implies $\t_1 = a \in \Q$, which cannot happen). Then $\tau_2 \in \spn_\Q(1, \tau_3)$. If $d=0$, then $\t_2 \in \Q$, otherwise $\t_2$ depends non trivially on $\t_3.$

A similar argument with $\begin{pmatrix}
0 & 0 \\ 1 & 0
\end{pmatrix} $ yields equations 
\begin{align*}
0 = a' + c'\tau_1 + d' \tau_2,\\
\tau_1 = b' + c'\tau_2 + d'\tau_3,
\end{align*}
where from the first equation $d' \not= 0$ and $\tau_2\in \spn_\Q(1, \tau_1)$. If $c'=0$, it follows that $\t_2 \in \Q$ and the second equation gives a non trivial linear relation between $\t_1$ and $\t_3$ (as $d' \not= 0)$, and we have proven our result. If $\t_2 \in \Q$, as the first case above, this argument still aplies, and we arrive at the desired conclusion.

Finally, in the more general case we have $\t_2 \in \spn_\Q(1, \t_1)\cap \spn_\Q(1, \t_3)$ with non trivial dependence on $\t_1$ and $\t_3$. Solving for $\tau_2$ gives a non zero linear equation that involves $1, \tau_1, \tau_3$. The conclusion follows.
\end{proof}

\subsection{End of proof of Theorem \ref{TeoBigTeo}}
We finish the proof of Theorem \ref{TeoBigTeo}. If $A_{\bt} \simeq E^2$, with $E$ a CM elliptic curve, as we discuss at the beginning of Section \ref{sec:SecRankHigh}, we have
\[
\trqjq = 1.
\]
By Theorem \ref{TeoShiodaMitani}, $A_{\bt} \in K(N)$ for some $N$, and hence by \cite[Theorem 12]{KaniThemodulispace}, it belongs to one of the curves in the collection described in Lemma \ref{LemaFamilyModCurves}, then generically in this curve $\rank \L_{\bt} = \Llin_{\bt} = 2$, and hence, if $A_{\tilde{\bt}} \in C$ is defined over $\Qbar$, 
\[
  \min_{\Sp_4(\Z)} \trqjqtil \leq 1.
\]

\appendix

\section{On a "simultaneous Humbert's lemma" for QM abelian surfaces}
Recall the collection of modular curves from Section \ref{sec:SecFamily}. For $A_{\bt}$ belonging to one of those, we have explicitly given an element in the $\Sp_4(\Z)$-orbit such that $\rank \Llin_{\bt} = 2$. We can rephrase that as an answer to a different question. Humbert's lemma Proposition \ref{PropHumbertsLemma} states that if $\bt \in \H_2$ solves a HSR, there is always an element in the orbit solving a linear HSR. For the elements in the locus in $\A_2$ covered by the collection in Section \ref{sec:SecFamily}, they solve \emph{two} linear ones. Then the question that we ask is: given $\bt \in \H_2$ solving two Humbert singular relations, can we always find an element $M \in \Sp_4(\Z)$ such that $\bt' = M\bt$ solves two linear HSR \emph{simultaneously?} We made the conjecture in Remark \ref{RemConjNonSimHumbForShimCurves} that it was not possible for abelian surfaces with $\End_0(A)$ a division quaternion algebra.

The explicit embedding from Proposition \ref{PropHas95} satisfies only $\rank \Llin_{\bt} =1$. We claim that, under some extra conditions, $\rank \Llin_{\bt} =1$ holds for the rest of the $\Sp_4(\Z)$-orbit. In other words, there is no "simultaneous Humbert lemma" for Shimura curves that allow us to linearize two HSR. We will now give a proof.

Our strategy is to translate the set-up to a Hilbert modular surface, as both Shimura and modular curves are compatible with real multiplication by the ring of integers of some quadratic field. We prove that both types of curves analytically correspond to Hirzebruch-Zagier divisors in $\H_1^2$. Furthermore, linear HSR will correspond to linear Hirzebruch-Zagier divisors. In this context, we have at our disposition a criterion that allows us to distinguish both types of curves by the coefficients or the divisors, and in particular implies that for compact curves (Shimura curves) the divisor cannot be linear. 

We set up the notation for Hilbert modular surfaces from \cite{Milio-Robert}, and the convention for Hirzebruch-Zagier divisors from \cite[Chapter V]{vanHMS}. 

Set a fundamental discriminant $\Delta > 0$ with $\Delta \equiv 1 ($ mod $4$) and consider the real quadratic field $K = \Q(\sqrt{\Delta})$. There is an explicit embedding from the (symmetric) Hilbert modular space to the Humbert surface $\Hum_{\Delta}.$

Consider the following notations.
\begin{itemize}
\item $\O_{K}$ for the ring of integers of $K$;
\item $w$ for the generator of $\O_{K}$, \textit{i.e.} $w = \frac{1 + \sqrt{\Delta}}{2}$ and $\{1, w\}$ is a $\Z$-basis of $\O_K$;
\item $\partial_{K} = \sqrt{\Delta} \O_{K}$ and $\partial_{K}^{-1} =\frac{1}{\sqrt{\Delta}} \O_{K}$ for the different and codifferent of $K$. Remark that the norm of the different ideal $\No(\partial_K) = \Delta$;
\item for $\al \in K$, write $\bar{\al}$ for the non-trivial Galois conjugate;
\item for $\mathfrak{a}$ a fractional ideal in $K$,
\begin{align*}
\Gamma( \mathfrak{a}) & := \SL_2(\O_K \oplus \mathfrak{a}) = 
\left\{ \begin{pmatrix}
a & b \\ c & d
\end{pmatrix} \in \SL_2(K) \; a,d \in \O_{K}, \: b \in \mathfrak{a}^{-1}, c \in \mathfrak{a} \right\} = \\
&= \begin{pmatrix}
\O_K & \mathfrak{a}^{-1} \\ \mathfrak{a} & \O_K
\end{pmatrix} \cap \SL_2(K),
\end{align*}
and we fix $\mathfrak{a} = \partial_K$;
\item the variables of the Hilbert space are denoted $\bz = (z_1, z_2) \in \H_1^2$.
\end{itemize}

Consider the matrix $R = \begin{pmatrix}
1 & w \\1 & \bar{w}
\end{pmatrix}$ and the map
\begin{align*}
\phi: \H_1^2 &\to \H_2 \\
\bz & \mapsto {}^t\!R \begin{pmatrix}
z_1 & 0 \\0 & z_2
\end{pmatrix} R = \begin{pmatrix}
z_1 + z_2 & z_1w + z_2\bar{w} \\ z_1w + z_2\bar{w}  & z_1w^2 + z_2\bar{w}^2
\end{pmatrix}.
\end{align*}
Then by \cite[Proposition 2.11]{Milio-Robert}), it defines an embedding between the symmetric Hilbert modular surface associated to $K$ and $\Hum_{\Delta}$ (after properly setting an embedding $\Gamma = \Gamma(\partial_K) \to \Sp_4(\Z)$) and it is independent of the choice of $\Z$-basis for $\O_K$). Setting $\bt = \phi(\bz)$, it satisfies the normalized HSR (note the change of convention in \cite{Milio-Robert})
\[
\frac{\Delta - 1}{4} \tau_1 + \tau_2 - \tau_3 = 0.
\]
We invert $\phi$ in $\Hum_{\Delta}$:
\begin{align*}
z_1 &= \frac{-\bar{w}}{\sqrt{\Delta}} \tau_1 + \frac{1}{\sqrt{\Delta}}\tau_2 =  \frac{-\bar{w}\tau_1 + \tau_2}{\sqrt{\Delta}},\\ 
z_2 &= \frac{w}{\sqrt{\Delta}}\tau_1 - \frac{1}{\sqrt{\Delta}}\tau_2 = \frac{w\tau_1 + \tau_2}{\sqrt{\Delta}}.
\end{align*}
Let us consider now $\bt$ with $A_{\bt}$ belonging to a Shimura curve in $\Hum_{\Delta}$. After applying Humbert's lemma to the relation of discriminant $\Delta$, can consider an intersection
\[
\begin{cases}
\frac{\Delta-1}{4} \tau_1 + \tau_2 - \tau_3 = 0\\
a'\tau_1 + b'\tau_2 + c\tau_3 + d(\tau_2 - \tau_1\tau_3) + e  = 0 \\ 
\end{cases}\implies  \begin{cases}
\frac{\Delta-1}{4} \tau_1 + \tau_2 = \tau_3\\
a\tau_1 + b\tau_2 - d\det(\bt) + e  = 0, \\
\end{cases}
\]
where $a = a' + c(\Delta-1)/4$, $b = b'+ c.$ 

\begin{lema} Under the map $\phi$, there is a bijection:
\begin{align*}
\left\{ \begin{array}{c}
(a,b,d,e) \in \Z^4: \\
 a\tau_1 + b\tau_2 + d\det\bt + e  = 0 
 \end{array}
 \right\} &\to \left\{ \begin{array}{c}
 (p,q, \gamma) \in \Z^2 \times \partial_K^{-1}: \\
  p\sqrt{\Delta} z_1z_2 - \bar{\gamma}z_1 + \gamma z_2 + \frac{q}{\sqrt{\Delta}} = 0
  \end{array}
   \right\}\\
(a,b,d,e) & \mapsto \left(d, e, \frac{1}{\sqrt{\Delta}}\left( a + b \bar{w} \right) \right).
\end{align*}
In particular, linear equations in $\bt$ correspond to linear equations in $\bz$. \label{LemmaLinEqtoLinEq}
\end{lema}
\begin{proof}
First, by definition of $\phi$ we have the following relation between determinants:
\[
\det\bt = z_1 z_2 (w - \bar{w})^2 = \Delta z_1 z_2.
\]

In one direction, setting $\tau_1 = z_1 + z_2$, $\tau_2 = wz_1 + \bar{w}z_2$:
\begin{align*}
a \tau_1 + b \tau_2 + d \det\bt + e &= 0\\
a\left( z_1 + z_2 \right) + b \left( wz_1 + \bar{w}z_2 \right) + d \Delta z_1 z_2 + e &= 0 \\
d \Delta z_1 z_2 + (a + bw)z_1 + (a + b\bar{w})z_2 + e &= 0\\
d \sqrt{\Delta}z_1 z_2 + \left( \frac{a + bw}{\sqrt{\Delta}} \right) z_1 + \underbrace{\left( \frac{a + b\bar{w}}{\sqrt{\Delta}} \right)}_{=: \gamma} z_2 + \frac{e}{\sqrt{\Delta}} &= 0
\end{align*}
Clearly $\gamma \in \partial_{K}^{-1}$, and $-\bar{\gamma} = -\frac{a + bw}{-\sqrt{\Delta}} = \frac{a + bw}{\sqrt{\Delta}}$. 

In the other direction,
\begin{align*}
p\sqrt{\Delta} z_1z_2 - \bar{\gamma}z_1 + \gamma z_2 + \frac{q}{\sqrt{\Delta}} &= 0\\
p\sqrt{\Delta} \frac{\det\bt}{\Delta} - \bar{\gamma}\left(\frac{-\bar{w}\tau_1 + \tau_2}{\sqrt{\Delta}} \right) + \gamma \left(  \frac{w\tau_1 + \tau_2}{\sqrt{\Delta}} \right) + \frac{q}{\sqrt{\Delta}} &= 0\\
\left(\gamma w + \overline{\gamma w} \right)\tau_1 -\left(\bar{\gamma} + \gamma \right)\tau_2 + p\det(\bt) + q &= 0,
\end{align*}
and we just need to check that the coefficients are integers (remark that neither $\gamma$ nor $w\gamma$ are algebraic integers). Write $\gamma = (m + nw)/\sqrt{\Delta}$, then:
\[
\gamma + \bar{\gamma} = \frac{m + nw -m - n\bar{w}}{\sqrt{\Delta}} = \frac{n(w - \bar{w})}{\sqrt{\Delta}} = n \in \Z,
\]
and the same argument applies to $w\gamma$, as it also belongs to $\partial_K^{-1}$.
\end{proof}

Let us arrange the coefficients of the equations in $\bz$ in a matrix
\[
B = \begin{pmatrix}
p\sqrt{\Delta} & \gamma \\ -\bar{\gamma} & \frac{q}{\No(\partial_K)} \sqrt{\Delta}
\end{pmatrix},
\]
with $p,q \in \Z, \gamma \in \partial_K^{-1}$. Then this matrix is what in \cite[Chapter V, Definition(1.2)]{vanHMS} is called a \emph{skew-hermitian matrix integral with respect to} $\partial_K$. Its determinant satisfies:
\[
pq + \No(\gamma) \in \frac{1}{\No(\partial_K)} \Z.
\]
One can think of $\Delta \det (B)$ as an analog to the discriminant of the HSR, in the sense that the action of $\Gamma$ changes the equation but respects the determinant. Following \cite{vanHMS}, we define the subset of $\H_1^2:$
\[
T(M) = \bigcup_{\substack{ B = \begin{pmatrix}
p\sqrt{\Delta} & \gamma \\ -\bar{\gamma} & \frac{q}{\No(\partial_K)} \sqrt{\Delta}
\end{pmatrix} \\ \det(B) = \frac{M}{\No(\partial_K)}  }}    \{ (z_1, z_2) \in \H_1^2: \; p\sqrt{\Delta} z_1z_2 - \bar{\gamma}z_1 + \gamma z_2 + \frac{q}{\No(\partial_K)} \sqrt{\Delta} = 0 \},
\]
which is $\Gamma$-invariant and defines a divisor on $\H_1^2$, the \emph{Hirzebruch-Zagier divisor of discriminant} $M$.  We remark that the standard definition of these divisors requires $\det B$ positive, and the ones considered in Lemma \ref{LemmaLinEqtoLinEq} do not necessarily satisfy that. However, the relevant results to us (namely the discussion between \cite[Chapter V, Lemma 1.4 and Proposition 1.5 ]{vanHMS}) do not require the determinant to be positive.

As we have said, this divisor can define a Shimura or modular curve, and one can tell them apart because the former are naturally compact (equivalently, they do not meet the resolution of the cusps of the compactification of the Hilbert surface). More precisely, it follows from \cite[Chapter V, Proposition 1.5]{vanHMS} (and from the prior discussion, which does not require $\det(B)>0$), that the indefinite quaternion algebra over $\Q$ corresponding to $T(M)$ is
\[
\left( \frac{\Delta, \frac{-M}{\No(\partial_K)\Delta }}{\Q} \right). 
\]

We are finally able to prove that Shimura curves do not admit a simultaneous Humbert lemma. Arguing by contradiction, assume that it admits a description by intersection of two linear HSR. 

\begin{vacio}
Assume that one of them is a normalized HSR with respect to a fundamental discriminant $\Delta \equiv 1$ (mod $4$). \label{vacio} 
\end{vacio}

By Lemma \ref{LemmaLinEqtoLinEq}, the linear HSR gives a linear equation for the Hirzebruch-Zagier divisor corresponding to the Shimura curve in the Hilbert modular surface of discriminant $\Delta$. The corresponding matrix is of the shape
\[
B = \begin{pmatrix}
0 & \gamma \\ - \bar{\gamma} & \frac{q}{\sqrt{\Delta}}
\end{pmatrix},
\]
with $\det(B) = \gamma \bar{\gamma} \not= 0.$ Consider the quaternion algebra
\[
\left(\frac{\Delta, \frac{-\det(B)}{\Delta }}{\Q} \right) = \left( \frac{\Delta, \frac{-\gamma \bar{\gamma}}{\Delta}}{\Q}\right) = \left( \frac{\Delta, \No(\al)}{\Q} \right) 
\]
where $\al := \gamma /\sqrt{\Delta}$ (remark the change of sign, because $\No(\sqrt{\Delta}) = -\Delta$). It is finally enough to see that the quaternion algebra splits over $\Q$ to arrive to a contradiction. 

\begin{lema} Consider $\Delta >0$ a fundamental discriminant, $K = \Q(\sqrt{\Delta})$ and $\al \in K^{\times}$. It follows
\[
\left( \frac{\Delta, \No(\al)}{\Q} \right) \cong \M_2(\Q)
\]
\end{lema}
\begin{proof}
Consider the standard $\Q$-basis $\{1,I,J,IJ\}$ with $I^2 = \Delta$, $J^2 = \No(\al)$ and $IJ = -JI$. We can then identify $I=\sqrt{\Delta}$. By standard theory of quaternion algebras over $\Q$, the statement is equivalent to the quaternion algebra having a non-trivial element of norm zero. We can explicitly find one: write $\al = a + b \sqrt{\Delta}$ with $a,b \in \Q$. Then 
\[
\No(\al) = a^2 -\Delta b^2.
\]
Consider the element $\mu := a + bI + J$. Then
\[
(a + bI +J)(a - bI - J) = a^2 - (bI)^2 - J^2 = a^2 - \Delta b^2 - \No(\al) = 0.
\]
\end{proof}

The condition $\Delta \equiv 1$ (mod $4$) in Assumption \ref{vacio} could be removed by repeating this analysis for the quadratic fields of discriminant divisible by four. The true impositions are that one of the linear HSR is \emph{normalized} and with respect to a \emph{fundamental} discriminant. Hence this proof does not rule out that the linear HSR could be not normalized, or normalized with respect to an order in a quadratic field that it is not the ring of integers.

\printbibliography

\end{document}